\documentclass[]{ijmart}

\usepackage[utf8]{inputenc}
\usepackage[T1]{fontenc}

\usepackage[english]{babel}
\usepackage[colorlinks]{hyperref}

\usepackage[]{caption}
\title[Arithmetical Hierarchy of Stability]{Arithmetical Hierarchy of the Besicovitch-Stability of Noisy Tilings}

\author[Léo Gayral]{Léo Gayral\thanks{
This work is an extension upon preliminary results
first introduced as an exploratory paper by the first author at Automata 2021~\cite{Gay21}.
}}
\email{\href{mailto:leo.gayral@math.cnrs.fr}{leo.gayral@math.cnrs.fr}}
\urladdr{\href{https://lgayral.pages.math.cnrs.fr/en/}{lgayral.pages.math.cnrs.fr/en}}

\author{Mathieu Sablik}
\email{\href{mailto:msablik@math.univ-toulouse.fr}{msablik@math.univ-toulouse.fr}}
\urladdr{\href{https://www.math.univ-toulouse.fr/\~msablik/index.html}{math.univ-toulouse.fr/\raisebox{0.5ex}{\texttildelow}msablik}}

\address{\hspace{0pt}\\
Université Toulouse III - Paul Sabatier\\
Institut de Mathématiques de Toulouse\\
118, route de Narbonne\\
F-31062 Toulouse Cedex 9\\
France
}

\date{} 
\issueinfo{VOL}{NUM}{MONTH}{2023}
\doiinfo{10.1007/DOI-NUMBER}

\newcommand{\A}{\mathcal{A}}
\newcommand{\B}{\mathcal{B}}
\newcommand{\F}{\mathcal{F}}
\newcommand{\M}{\mathcal{M}}
\newcommand{\W}{\mathcal{W}}
\renewcommand{\epsilon}{\varepsilon}
\renewcommand{\phi}{\varphi}

\usepackage{dsfont}
\newcommand{\N}{\mathds{N}}
\newcommand{\Z}{\mathds{Z}}
\newcommand{\Q}{\mathds{Q}}
\newcommand{\Prob}{\mathds{P}}
\newcommand{\R}{\mathds{R}}

\newcommand{\meas}[1]{\widehat{\delta_{#1}}}
\newcommand{\Ball}{\overline{\mathrm{B}}}

\usepackage{stmaryrd}
\usepackage{amsmath}
\usepackage{amssymb}

\newtheorem{theorem}{Theorem}[section]
\newtheorem{corollary}[theorem]{Corollary}
\newtheorem{lemma}[theorem]{Lemma}
\newtheorem{proposition}[theorem]{Proposition}

\theoremstyle{definition}
\newtheorem{definition}{Definition}[section]
\newtheorem{remark}{Remark}[section]

\usepackage{xcolor}

\usepackage{float}
\usepackage{graphicx}
\usepackage[font=small]{caption}

\newcommand{\figaligntop}{\vspace{0pt}}
\newcommand{\figalignmid}{\vspace{-6pt}}
\newcommand{\figalignbot}{\vspace{0pt}}

\fancypagestyle{firstpage}{
  \fancyhf{}
  \chead{\tiny
   \hspace{1pt}\\[0.5ex] 
   \hspace{1pt}}
   \cfoot{\thepage}}

\begin{document}

\begin{abstract}
The purpose of this article is to study the algorithmic 
complexity of the Besicovitch stability of noisy
subshifts of finite type,
a notion studied in a previous article~\cite{GaySa21}.
First, we exhibit an unstable aperiodic tiling,
and then see how it can serve as a building block
to implement several reductions
from classical undecidable problems on Turing machines.
It will follow that the question of stability of
subshifts of finite type is undecidable,
and the strongest lower bound we obtain
in the arithmetical hierarchy is $\Pi_2$-hardness.
Lastly, we prove that this decision problem,
which requires to quantify over
an uncountable set of probability measures,
has a $\Pi_4$ upper bound.
\end{abstract}

\maketitle
\tableofcontents


\section{Introduction}

Let $\A$ a finite alphabet.
A subshift of finite type (\emph{SFT}), denoted $\Omega_\F$,
is a set of $\A$-colourings of $\Z^d$ induced by
a finite set of forbidden patterns $\F$ which
cannot appear in any configuration.
One of the main topics of interest in the study
of multidimensional SFTs is how a global structure
can emerge from local rules.
In particular, \emph{aperiodic} SFTs have been studied
by Berger~\cite{Ber66}, Robinson~\cite{Rob71}
and Kari~\cite{Kari96} among others.
One of the most useful properties
of the Robinson tiling
is that its hierarchical structure leaves room for
a \emph{relatively} easy embedding of
Turing machines into it~\cite{Rob71,JeanVa20}.

In the last decade, a lot of studies focused on
the links between dynamical properties of SFTs
and their algorithmic complexity.
The values taken by some dynamical
invariants can be characterised as 
some classes of (non-)computable values:
possible entropies~\cite{Hochman-Meyerovitch-2010},
or dimension entropies~\cite{Meyerovitch-2011},
subactions~\cite{AuSab10,DuRoShe12,Hoch09},
possible periods~\cite{Jeandel-Vanier-2014},
or some classes of SFTs~\cite{Westrick-2017}\dots{}
These works help to understand the limits of
what global behaviours can be enforced by local rules.

These classes of \emph{numbers} relate
to the arithmetical hierarchy of computable \emph{sets}
through the identification between $x\in\R$
and the interval $\{q\in\Q,q<x\}$.
Another way to highlight the complexity of tilings
is then to understand the complexity of
a decision problem about a dynamical property of the SFTs.
These problems are usually undecidable,
but may fit into the arithmetical (or analytical) hierarchy.
Regarding the arithmetical hierarchy, 
the Domino problem is $\Pi_1$-complete~\cite{Rob71},
the conjugacy problem is $\Sigma_1$-complete
and the factorisation problem is $\Sigma_3$-complete
~\cite{Jeandel-Vanier-2015}\dots{}
Regarding the analytical hierarchy,
deciding whether a tiling has
a completely positive topological entropy or not
is $\Pi_1^1$-complete~\cite{Westrick-2022},
in dimension $d\geq 4$ the aperiodic Domino problem 
is $\Pi_1^1$-complete~\cite{Hellouin-2022}\dots{}
To obtain these results,
the proofs always involve the embedding
of Turing machines into \emph{complex}
(and aperiodic) tilings.
This is interesting since few natural problems
(not directly related to a computation model)
are known to be complete in these hierarchies.

In this article, we study the algorithmic complexity
of the Besicovitch-stability of noisy SFTs.
In a previous article~\cite{GaySa21},
we introduced this notion of stability
using the Besicovitch distance,
which quantifies the closeness between measures
through the average frequency of differences
between their configurations.
This framework is a natural bridge from
the notion of stability described by
Durand, Romashchenko and Shen~\cite{DuRoShe12}
to ergodic theory,
with a viewpoint focusing more on measure theory.
The purpose is to understand if SFTs are stable
in the presence of noise,
if computations can survive if a small
proportion of forbidden patterns is permitted.
Such studies already exist
for cellular automata~\cite{Ga01}
or Turing machines~\cite{AsaCol05}.
A digest of this framework will be introduced
in Section~\ref{sec:Framework},
followed by a few notions about undecidability
and the arithmetical hierarchy.

In the aforementioned article~\cite{GaySa21},
we proved a simple computable criterion
(using a word automaton)
to decide stability for one-dimensional SFTs.
Then, we proved the existence of both stable
\emph{and} unstable SFTs in \emph{any} dimension, 
and a specific variant of the Robinson tiling
was proven to be stable;
before this, the only known stable aperiodic tilings
were complex constructions
that can be repaired locally,
which is not the case for this variant of the
Robinson tiling~\cite{BaDuJean10,DuRoShe12,Taati}.
However, the interface between stable and unstable
examples in general was yet to be seen.

In this article, we will prove that a known
two-coloured Robinson tiling is unstable in
Section~\ref{sec:RedBlackRobinson},
and describe a general framework to obtain stability
for some quasi-periodic SFTs
in Section~\ref{sec:AperiodicStability}.
By iterating upon both the stable and unstable constructions,
we will step-by-step craft simulating tilings to show
that deciding if a SFT is stable is $\Pi_1$-hard,
$\Sigma_1$-hard and finally $\Pi_2$-hard
in Section~\ref{sec:UndecidableStability}.

After this, we will obtain a $\Pi_4$ upper bound
for stability in Section~\ref{sec:UpperBound}. 
This bound may be surprising \emph{a priori}
since the definition of stability requires to
quantify over uncountable sets
(of translational-invariant probability measures).
To obtain such a bound we will dig deeper into
the technicalities of computable analysis on measures,
to rewrite the stability property
using only elements from a countable basis.
This section is independent of the previous constructs
for the lower bounds,
and relies only on the definitions
of Section~\ref{sec:Framework}.  


\section{General Framework} \label{sec:Framework}

In this section, we define the general framework for the rest of the paper.

First, we introduce noisy SFTs and stability,
which were defined more in-depth
in a previous paper~\cite[Sections 2 and 3.1]{GaySa21}.
This subsection explains most of the notations used later on,
and provides a baseline of ergodic theory for readers with a computer science background
in particular.

Second, we define what decidability and the arithmetical hierarchy mean in our context,
so that readers with a mathematical background in particular can still follow the rest.

\subsection{Noisy SFTs and Besicovitch Stability} \label{sec:Besicovitch}

\begin{definition}[Subshift of Finite Type]
Let $\A$ be a finite alphabet, and denote $\Omega_\A:=\A^{\Z^d}$,
endowed with the product topology and corresponding Borel algebra.
Let $\F$ be a finite set of \emph{forbidden} patterns $w\in \A^{I(w)}$,
defined on finite windows $I(w)\Subset \Z^d$.
A SFT is the set $\Omega_\F$ induced by $\F$ as follows:
\[
\Omega_\F:=\left\{ \omega\in \Omega_\A ,\forall w \in\F,
\forall k\in\Z^d,\sigma_k(\omega)|_{I(w)}\neq w \right\} ,
\]
\textit{i.e.}\ configurations of the SFT are such that
no forbidden pattern occurs.

This set is $\sigma$-invariant,
invariant
for any translation $\sigma_k$ (with $k\in\Z^d$),
defined as $\sigma_k:\left(\omega_l\right)_{l\in\Z^d} 
\mapsto \left(\omega_{k+l}\right)_{l\in\Z^d}$.
Thus, if we denote $\left(e_i\right)_{1\leq i \leq d}$
the canonical basis of $\Z^d$,
$\left(\Omega_\F,\sigma_{e_1},\dots,\sigma_{e_d}\right)$
is a commutative dynamical system.
\end{definition}

Now, we twist this notion to include noise through \emph{obscured} cells:

\begin{definition}[Noisy SFT]
Consider the alphabet $\widetilde{\A}=\A\times \left\{0,1\right\}$,
with the identification $\A\approx\A\times \left\{0\right\}$.
Formally, we denote $\pi_1:\widetilde{\A}\to \A$
and $\pi_2:\widetilde{\A}\to\left\{0,1\right\}$ the canonical projections.
We can likewise define the set of forbidden patterns $\widetilde{\F}
:=\left\{ \left(w, 0^{I(w)}\right)\in\widetilde{A}^{I(w)}, w\in\F\right\}$
and the corresponding SFT $\Omega_{\widetilde{\F}}$
on $\widetilde{A}$.
\end{definition}

In general, if $\mu$ is a measure on $\Omega$
and $\phi:\Omega\to\Omega'$ is a measurable mapping,
we can define the pushforward measure $\phi^*(\mu)$ on $\Omega'$,
such that for any measurable set $A\subset \Omega'$,
we have $\left[\phi^*(\mu)\right](A)=\mu\left(\phi^{-1}(A)\right)$.

\begin{definition}[Noisy Probability Measures]
A measure $\mu$ is $\sigma$-invariant if for any $k\in\Z^d$,
the pushforward measure $\sigma_k^*(\mu)$ is equal to $\mu$.
Denote $\M_\F$ the set of $\sigma$-invariant probability measures
supported by $\Omega_\F$.

Let $\B:=\left\{ \B(\epsilon)^{\otimes \Z^d}, 0\leq \epsilon\leq 1\right\}$
be the class of Bernoulli noises.
Define:
\[
\widetilde{\M_\F^\B}(\epsilon) :=\left\{
\lambda \in \M_ {\widetilde{\F}},
\pi_2^*(\lambda)\in\mathcal{B}\textrm{ and }\pi_2^*(\lambda)([1])\leq \epsilon\right\} .
\]
Likewise, $\M_\F^\B(\epsilon) :=
\pi_1^*\left(\widetilde{M_\F^\B}(\epsilon) \right)$
consists of probability measures on $\Omega_\A$.
\end{definition}

The measures of $\widetilde{\M_\F^\B(\epsilon)}$
have a low probability of containing
obscured cells in a given finite window.
However, we still need a way to globally quantify the structural effect of these few local errors:

\begin{definition}[Besicovitch Distance]
We define $d_H$ the Hamming-Besicovitch pseudo-distance
on $\Omega_\A$ as $d_H(x,y)=\varlimsup\limits_{n\to\infty}
d_n\left(x|_{B_n},y|_{B_n}\right)$, with the Hamming distances
$d_n(u,v)=\frac{1}{(2n+1)^d}\#\left\{k\in B_n,
u_k\neq v_k\right\}$ and $B_n:=\llbracket -n,n\rrbracket^d$.

A \emph{coupling} (or \emph{joining})
between two measures $\mu$ on $\Omega_{\A_1}$ and $\nu$ on $\Omega_{\A_2}$ 
is a measure $\lambda$ on $\Omega_{\A_1\times \A_2}$
such that $\pi_1^*(\lambda)=\mu$ and $\pi_2^*(\lambda)=\nu$.
Denote $J(\mu,\nu)$ the set of such couplings, and more generally
$J(U,V)=\bigcup_{\mu\in U,\,\nu\in V} J(\mu,\nu)$.
The Besicovitch distance between two $\sigma$-invariant measures is then:
\[
d_B(\mu,\nu):=\inf\limits_{\lambda\in J(\mu,\nu)}\int d_H(x,y)\mathrm{d}\lambda(x,y) .
\]
\end{definition}

By $\sigma$-invariance of the measure $\lambda$,
an ergodic theorem~\cite[Chapter 6]{Kre11} gives us a link between global and local scales through
$\int d_H(x,y)\mathrm{d}\lambda(x,y)=\lambda\left(\left[x_0\neq y_0\right]\right)$
with the cylinder set $\left[x_0\neq y_0\right]:=\left\{ (x,y)\in\Omega_\A^2,x_0\neq y_0\right\}$.
This equivalent definition of the distance can be in particular found as the distance $\overline{d}$
in \emph{Ergodic Theory via Joinings}~\cite[Chapter 15]{Gla03}.

For two ergodic measures, $d_B$ quantifies how well we can align
their generic configurations so that they coincide on a high density subset of $\Z^d$.
Using this distance, we can intuitively define stability as follows:

\begin{definition}[Stability]
The SFT $\Omega_\F$ is stable (for $d_B$ on $\B$)
if there is a non-decreasing $f:[0,1]\to \mathds{R}^+$,
continuous in $0$ with $f(0)=0$, such that:
\[
\forall \epsilon\in [0,1],\sup\limits_{\mu\in \M_\F^\mathcal{B}(\epsilon)}
d_B\left(\mu, \M_\F\right) \leq f(\epsilon) .
\]
\end{definition}

The general idea to keep in mind afterwards is that
this framework allows us to compare the average distance between configurations,
hence we will always go back to \emph{generic} configurations in some sense,
and compare these with $d_H$ to obtain a bound for $d_B$.

Now that stability has been defined,
we want to study its computational complexity.
As we will see later on, this problem is actually undecidable,
so we will want to see \emph{how much} undecidability it contains.
This is why we now need to introduce the notion of arithmetical hierarchy,
which allows for a classification of the complexity of undecidable problems.

\subsection{Decidability and the Arithmetical Hierarchy}

The goal of this subsection is to introduce the general vocabulary and key ideas,
so we will not plunge deep into the formalism, but we refer the interested reader
to the classical books by Rogers~\cite{Rog87} or Soare~\cite{Soa99}.
A less formal introduction on the topic can also be found on the
mathematical blog Rising Entropy~\cite{RisingEntropy20}.

A \emph{problem} is formally defined as
a subset of integers 
$P\subset\N$, usually described implicitly as the 
set of integers satisfying some
mathematical property.
Such a problem is said to be \emph{decidable} if there exists an algorithm
(or more formally a Turing machine) that answers in finite time when asked whether
$x\in\N$ belongs to $P$ or not.
If $P$ cannot be decided, it is called undecidable.

This notion (and the following ones) naturally extends to any countable space
that can be explicitly encoded into $\N$, such as $\Z^d$ for $d\geq 2$,
or the space of \emph{finite collections of (forbidden) patterns} $\F$.
Hence, we define $P_{stab}$ as the set of families of forbidden patterns $\F$
that induce a stable SFT.
The goal of the arithmetical hierarchy is to further classify these undecidable problems.

\begin{definition}[$\Pi_k$ and $\Sigma_k$ Problems]
We say that $P\in\Pi_k$ (with $k\in\N$) if we have $\phi\left(x,n_1,\dots,n_k\right)$
a computable algorithm on $\N^{k+1}$ such that $x\in P$
\emph{iff} the following formula holds true:
\[
\underset{k \text{ alternating quantifiers starting with }\forall}
{\underbrace{\forall n_1\in\N,\exists n_2\in \N,\forall n_3\in\N,\dots,}}\,
\phi\left(x,n_1,\dots,n_k\right) .
\]
Likewise, we say that $P$ is $\Sigma_k$ if we have the analogous property but starting
with an $\exists$ quantifier.
Note in particular how $\Pi_0=\Sigma_0$ simply describes decidable questions.
\end{definition}

It follows directly from the definition that
$\Pi_k\cup \Sigma_k \subset \Pi_{k+1}\cap\Sigma_{k+1}$,
and this inclusion is actually strict.

\begin{definition}[$\Pi_k$-hardness]
At last, we say that a problem $P$ is $\Pi_k$-hard if,
for any problem $Q\in\Pi_k$, there exists a computable reduction function $\phi:\N\to\N$ such that
$x\in Q$ iff $\phi(x)\in P$.
A problem $P$ is then $\Pi_k$-complete if $P\in\Pi_k$ and it is $\Pi_k$-hard.
\end{definition}

Notoriously, the \emph{halting} problem $P_{halt}$
(\textit{Does a Turing machine $M$ halt on the empty input?}) is $\Sigma_1$-complete,
and the \emph{totality} problem $P_{total}$ (\textit{Does $M$ halt on all of its inputs?})
is $\Pi_2$-complete~\cite[Part A, Chapter IV, Theorem 3.2]{Soa99}.
In Section~\ref{sec:UndecidableStability}, we will establish a computable reduction
from these problems to $P_{stab}$ to obtain a lower bound on its computational hardness.

As the definition of both stability \emph{and} instability
presuppose that $\Omega_\F\neq \emptyset$
(equivalent to the complementary of the halting problem, hence $\Pi_1$-complete),
we will include this property in the requirements for having $\F\in P_{stab}$.
With this unambiguous definition, in Section~\ref{sec:UpperBound},
we will prove a $\Pi_4$ upper bound on
the computational complexity of $P_{stab}$.


\section{The Red-Black Robinson Tiling is Unstable} \label{sec:RedBlackRobinson}

\begin{figure}[H]
\centering
\figaligntop
\includegraphics[scale=.3]{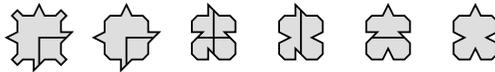}
\figalignmid
\caption{The $6$ basic Robinson tiles.\\
The leftmost one is called a bumpy-corner.}
\label{fig:VanillaRobinson}
\figalignbot
\end{figure}

Consider the Robinson tiling~\cite{Rob71} in Figure~\ref{fig:VanillaRobinson},
using the bumpy-corners variant (with diagonal interactions) instead of Wang tiles.
The tileset uses these $6$ tiles
and their rotations and symmetries, for a total of $32$ tiles in the alphabet.
The corresponding set of forbidden patterns is self-evident,
such that two laterally neighbouring tiles must have matching edges,
and each square of four tiles must use exactly one bumpy-corner to fill the hole in the middle.

\begin{figure}[H]
\centering
\figaligntop
\includegraphics[width=.38\linewidth]{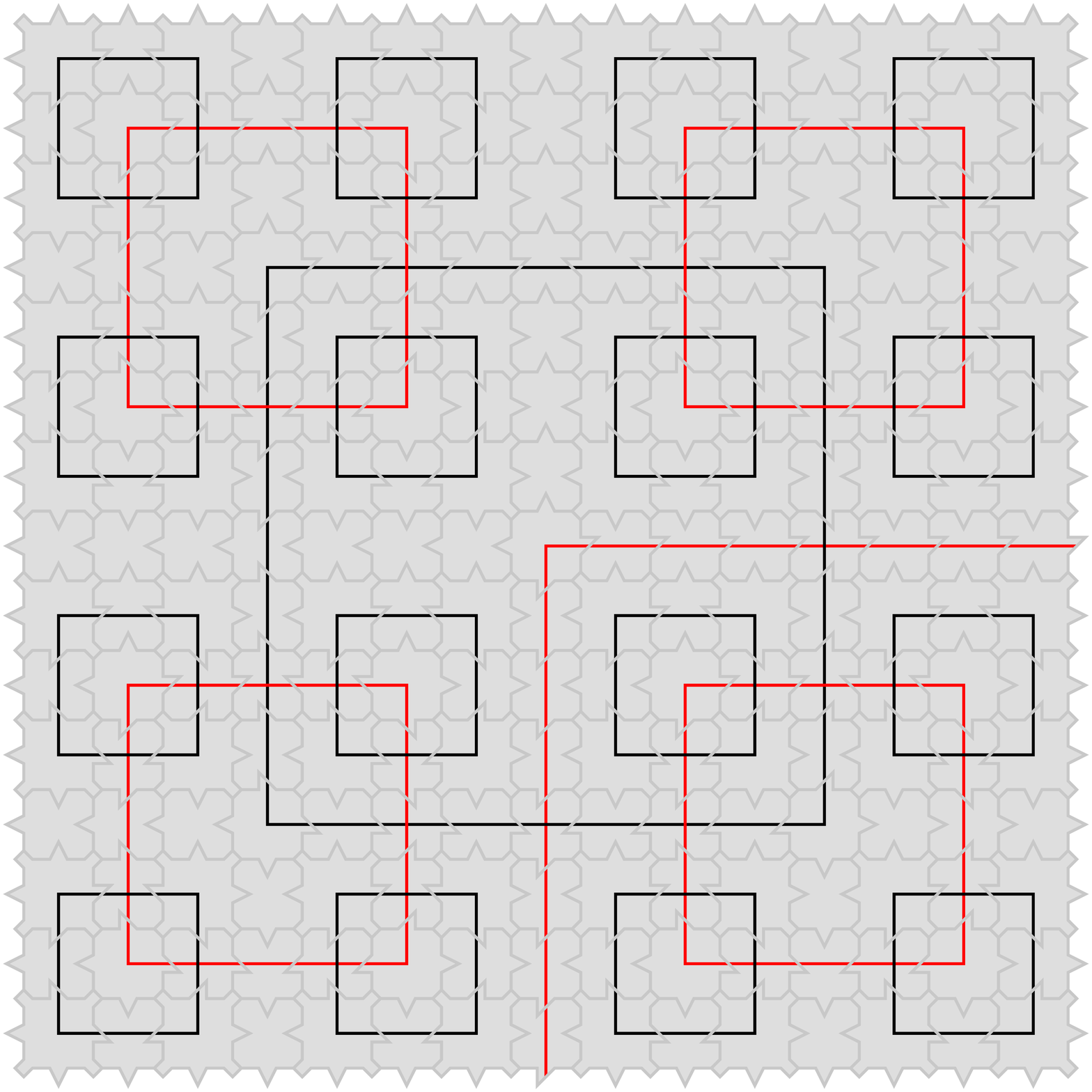}
\hspace{6mm}
\includegraphics[width=.38\linewidth]{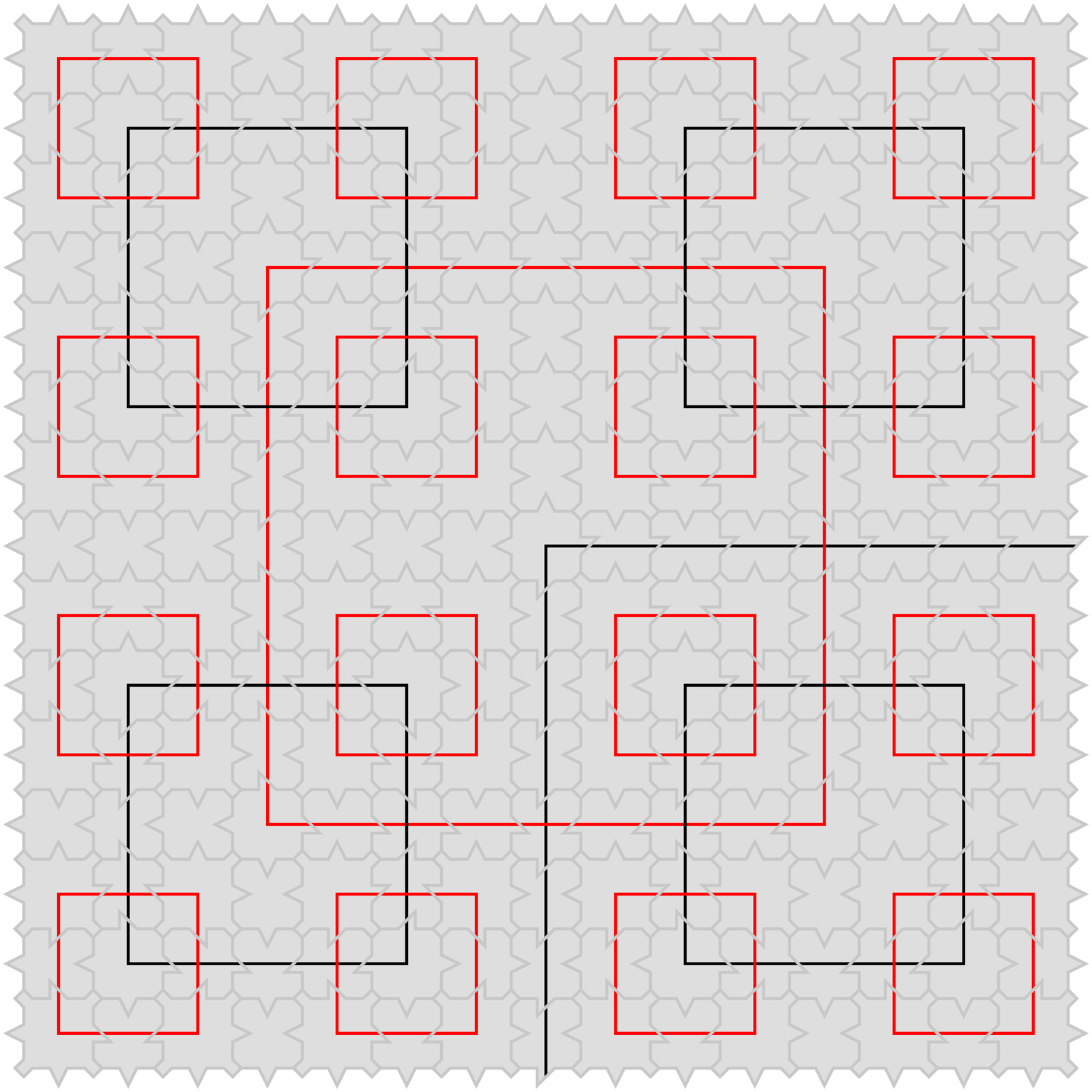}
\figalignmid
\caption{Alternating colours in Red-Black Robinson macro-tiles.}
\label{fig:RBStruct}
\figalignbot
\end{figure}

This tileset induces a self-similar hierarchical structure:
we first define the $1$-macro-tiles as the four rotated bumpy-corners tiles, 
and a $(n+1)$-macro-tile is then obtained
by sticking four $n$-macro-tiles in a square-like pattern,
around a central cross with two arms
(which itself has four possible orientations),
as in Figure~\ref{fig:RBStruct}.

In a previous paper~\cite[Theorem 7.9]{GaySa21}
we proved that an extension of this tileset,
enhanced to locally enforce the alignment of macro-tiles,
was stable with a polynomial speed $O\left(\sqrt[3]{\epsilon}\right)$.
Note that the Robinson tiling is not \emph{robust}
in the sense of Durand, Romashchenko and Shen~\cite{DuRoShe12},
so their anterior stability result
did not already apply to this tiling.

\begin{figure}[H]
\centering
\figaligntop
\includegraphics[scale=.3]{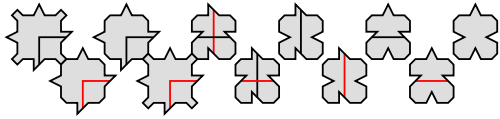}
\figalignmid
\caption{The $11$ basic Red-Black Robinson tiles.}
\label{fig:RBRobinson}
\figalignbot
\end{figure}

Here, we will use the two-coloured extension of this Robinson tileset
in Figure~\ref{fig:RBRobinson}, which naturally projects onto the previous tiling,
so all the structural properties of the Robinson tiling still hold,
and most notably aperiodicity.
We will denote $\A$ the tileset,
$RB$ the corresponding set of forbidden patterns,
and $\Omega_{RB}$ the resulting SFT.
Because $\A$ contains no tile with a monochromatic cross, only small crosses made of a
straight Red line crossing with a Black one,
any two squares of the same colour in the hierarchical structure of a tiling
do not intersect,
as we can see on the $5$-macro-tiles in Figure~\ref{fig:RBStruct}.
In Subsection~\ref{subsec:TuringRobinson},
these non-intersecting Red squares will be used to encode arbitrarily large
space-time diagrams of Turing machines.

For the rest of this paper,
a \emph{generic} Robinson tiling will refer to a configuration
without an infinite \emph{cut},
such that any two tiles of $\Z^d$ end up being in the same $n$-macro-tile
for big-enough values of $n$.
In such a \emph{generic} configuration $\omega\in\Omega_{RB}$,
by induction, the $n$-macro-tiles all have a central arm with the same colour.
In particular, a generic configuration will only contain
Red \emph{or} Black bumpy-corners, never both.

\begin{proposition} \label{prop:UnstableMeasure}
Let $\Omega_{RB}$ be the Red-Black Robinson tiling.
For any $\epsilon>0$,
there is $\mu\in\M_{RB}^\B(\epsilon)$ such that
$d_B\left(\mu,\M_{RB}\right) \geq \frac{1}{8}$.
Thus, the SFT is unstable.

\begin{proof}
The goal of this proof is to convert a generic tiling $\omega\in\Omega_{RB}$
into a random noisy tiling $\lambda_{\omega,b}$ on $\Omega_{\widetilde{RB}}$,
with $b$ a random variable on $\Omega_{\{0,1\}}$.
Using a generic Bernoulli noise $b$ in the input,
we will obtain a noisy tiling for which
its bumpy-corners are now half Red and half Black, which will yield the announced result since bumpy-corners have frequency $\frac{1}{4}$ in the Robinson tiling.

We will build this measure $\lambda$ iteratively,
as a limit of a locally-defined (thus trivially measurable) transformations.
At each step of the construction,
the actual monochromatic structure of the Robinson tiling will be preserved,
and only the colours will be mismatched,
so we may still consider $n$-macro-tiles in this structural sense,
even though they are not actually locally admissible.
We initialise $\lambda_1=\delta_{(\omega,b)}$ as a constant Dirac measure.

\begin{figure}[H]
\figaligntop
\centering
\includegraphics[width=.4\linewidth]{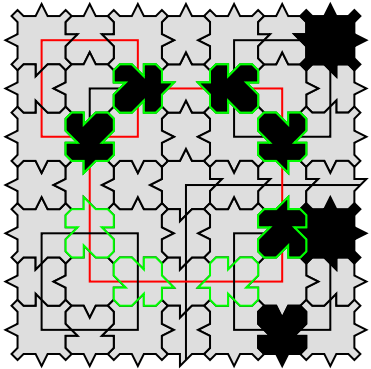}
\figalignmid
\caption{A locally admissible $3$-macro-tile with obscured cells.}
\label{fig:RBNoisy}
\figalignbot
\end{figure}

Let us now explain how we obtain $\lambda_2$ out of $\lambda_1$.
This transformation will be done independently
on \emph{each} of the $2$-macro-tiles of $\omega$.
We distinguish two cases, both illustrated in Figure~\ref{fig:RBNoisy},
where the black cells $c$ represent \emph{obscured} tiles with a noise $b_c=1$.
A macro-tile is said to be flippable if both of its bi-coloured crosses,
highlighted with green borders in the figure, are obscured tiles.
In such a situation, we will flip its colours (Black lines become Red and conversely)
with probability $\frac{1}{2}$, independently of the rest,
which still preserves the local rules inside the macro-tile.
In the figure, the top-left macro-tile is flipped,
the top-right macro-tile is flippable but not flipped,
and the two bottom macro-tiles are not flippable.

Likewise, we go from $\lambda_{n-1}$ to $\lambda_n$
by flipping independently at random any flippable $n$-macro-tile
(\emph{except} the two ends of the central arm that are ``after'' the bi-coloured crossed tiles,
which must match the colour of the yet-unflipped corresponding $(n+1)$-macro-tile).
This process guarantees that, if we denote $\omega'\sim \lambda_n$ the new colouring,
then $\left(\omega',b\right)\in \Omega_{\widetilde{RB}}$ almost-surely.

Notice how the highlighted cells that decide whether a given macro-tile is flippable
are disjoint for each macro-tile.
Hence, assuming that $b\sim\B(\epsilon)^{\Z^2}$ is a Bernoulli noise,
each macro-tile at each scale
is flippable with probability $\epsilon^2$, independently of the rest.
With such a choice of noise $b$, the weak-* limit $\lambda_{\omega,b}$ is well-defined.

Consider $G(\omega)\subset\Z^2$ the set of cells containing a bumpy corner in $\omega$.
For a given cell $c\in G$,
we denote by $\text{flip}_{c,n}$ the random variable equal to $1$
when the $n$-macro-tile containing $c$ is flippable.
Hence the variables $\text{flip}_{c,n}\sim\B\left(\epsilon^2\right)$ are \textit{iid}.
Conditionally to the event $\text{flip}_{c,n}=1$,
the colour of the cell $c$ is uniformly distributed in $\lambda_n$ after rank $n$.
Thence, by Borel-Cantelli lemma,
the colour of $c$ is uniformly distributed in $\lambda_{\omega,b}$.

Likewise, consider two distant cells $c,d\in G$.
As $d_\infty(c,d)\to\infty$,
the smallest rank $n_0(c,d)$ such that $c$ and $d$ belong to
the same $n$-macro-tile of $\omega$ goes to infinity.
The families $\left(\text{flip}_{c,n}\right)_{n<n_0}$
and $\left(\text{flip}_{d,n}\right)_{n<n_0}$ are independent,
and conditionally to the fact that both of these sequences contain at least a $1$,
the colours of cells $c$ and $d$ are independently uniform
(in the measures $\lambda_n$ after rank $n_0$, hence for $\lambda_{\omega,b}$).

Without loss of generality, assume $0\in G$, so that $G=(2\Z)^2$.
Then the family $\left(\text{colour of the cell }2c\right)_{c\in\Z^2}$
describes a $\sigma$-invariant ergodic dynamical system,
so that we may apply a pointwise ergodic theorem.
This implies that the \emph{frequency}
of both Black and Red bumpy-corners is generically equal
to $\frac{1}{2}$ in $\lambda_{\omega,b}$.
As $G$ has density $\frac{1}{4}$ in $\Z^2$,
we conclude that for almost-any $\omega'\sim\lambda_{\omega,b}$
and any generic $\omega_0\in\Omega_{RB}$ (with monochromatic bumpy-corners),
we have the bound
$d_H\left(\omega_0,\omega\right)\geq\frac{1}{2}\times\frac{1}{4}=\frac{1}{8}$
assuming bumpy-corners overlap between the two configurations,
and even a $\frac{1}{2}$ bound if they are misaligned.

We can conclude the proof by averaging $\lambda_{\omega,b}$
over $\omega\sim\mu_0\in\M_{RB}$ (chosen independently from $b$),
which gives us at last a $\sigma$-invariant measure $\mu\in\M_{RB}^\B(\epsilon)$
that satisfies $d_B\left(\mu,\M_{RB}\right)\geq\frac{1}{8}$.
\end{proof}
\end{proposition}

The result still holds with the very same proof if we replace
the bi-coloured Robinson tiling by a bi-coloured variant of the
structurally enhanced Robinson tiling from our previous paper~\cite{GaySa21}.

However, as the proof relies heavily on flipping the colours of bumpy-corners,
by keeping only \emph{one} of the two colours specifically for this tile,
we obtain a stable tiling again.
This will be useful later on, when we want to encode Turing machines
into Robinson (which requires this bi-coloured setting) in a stable way.
In such situations, stability will follow from the result of the next section.


\section{Generalising Aperiodic Stability} \label{sec:AperiodicStability}

In order to prove the stable cases later on,
we will state a direct generalisation of one of the main results
in our previous article~\cite[Proposition 7.8]{GaySa21}.
This proposition was proven in the specific context of the enhanced Robinson tiling,
but we will here reformulate the result in a general framework for quasi-periodic tilings
with a well-behaved reconstruction function,
so that it applies as a black box to
the tilesets described in the next section.
This section is here mostly for the sake of technical completeness,
and can be skipped to focus on the core of the paper to which we go back right after.

\begin{definition}[Almost Periodic SFT]
Let $\Omega_\F$ be a SFT on the alphabet $\A$, and consider $p\in\N^*$ and $\rho>0$.

We say that $\Omega_\F$ is $\rho$-almost $p$-periodic if
there is a $p$-periodic ``grid'' $G\subset \Z^d$
(invariant under translations in $(p\Z)^d$)
of density at most $\rho$,
such that any configuration restricted to a translation of $G^c$ is made periodic.
By this, we mean that for any $\omega\in\Omega_\F$,
there is a unique translation of $G$
(given by a non-necessarily unique
$k\in \llbracket 0,p-1\rrbracket^d$)
such that $\omega|_{G^c+k}$ is $p$-periodic.

In this case, assuming $\square \notin \A$,
we can define $\omega^\square$ by overwriting $\omega|_{G+k}$
by the blank symbol $\square$.
Thence, $\Omega_\F^\square=\left\{\omega^\square,\omega\in\Omega_\F\right\}$
is a finite $p$-periodic SFT.
\end{definition}

The non-uniqueness of $k$ comes from the fact that, for example,
we may want to consider $G$ a $\frac{p}{2}$-periodic grid instead,
with some more redundancy in its structure.

\begin{definition}[$C$-Reconstruction Function]
Consider $\Omega_\F$ a $\rho$-almost $p$-periodic SFT
and $G$ the associated grid.

The SFT has the $C$-reconstruction property if,
for any locally admissible tiling $\omega$ of $B_{\left\lceil\frac{p}{2}\right\rceil+C}$
there is a unique translation of $G$ such that
$\left[\omega|_{B_{\left\lceil\frac{p}{2}\right\rceil}\cap\left(G^c+k\right)}\right]^\square$
(obtained by filling $B_{\left\lceil\frac{p}{2}\right\rceil}\cap\left(G+k\right)$ with $\square$ symbols)
is globally admissible in $\Omega_\F^\square$
(thence $\omega|_{B_{\left\lceil\frac{p}{2}\right\rceil}\cap\left(G^c+k\right)}$ is
globally admissible in $\Omega_\F$).
What's more, the translation of $G$ depends only on what happens in any $p$-square included
in the central window $B_{\left\lceil\frac{p}{2}\right\rceil}$
(which is either a $(p+1)$ or $(p+2)$-square depending on the parity).
\end{definition}

As $\Omega_\F^\square$ is $p$-periodic,
there is a unique choice
$\omega^\square\in\Omega_\F^\square$
of configuration that will match the pattern
$\left[\omega|_{B_{\left\lceil\frac{p}{2}\right\rceil}\cap\left(G^c+k\right)}\right]^\square$.

\begin{proposition}[Besicovitch Bound] \label{prop:BesicovitchBound}
Consider a $\rho$-almost $p$-periodic SFT with $C$-reconstruction.
Then, for any $\epsilon>0$ and $\mu\in\M_\F^\B(\epsilon)$,
we have the bound $d_B\left(\mu,\M_\F\right)\leq 48\left(2\left(C+
\left\lceil\frac{p}{2}\right\rceil\right)+1\right)^d\epsilon+\rho$.

\begin{proof}
The proof is really similar to the source result~\cite[Proposition 7.8]{GaySa21},
so we will just give the general idea.

Consider $\lambda\in\widetilde{\M_\F^\B}(\epsilon)$ and
$(\omega,b)\in\Omega_{\widetilde{\F}}$ a $\lambda$-generic noisy configuration.
A percolation argument~\cite[Proposition 5.6]{GaySa21} tells us that, almost-surely,
we can forget about the $\left(C+\left\lceil\frac{p}{2}\right\rceil\right)$-neighbourhood
of \emph{obscured} cells (cells $c\in\Z^d$ with $b_c=1$)
and still have a unique connected component of \emph{clear} cells ($b_c=0$)
with density of at least $48\left(2\left(C+\left\lceil\frac{p}{2}\right\rceil\right)+1\right)^d\epsilon$.

Each clear cell $c$ of this connected component is the center of a clear window $I_c$
of diameter $2\left\lceil\frac{p}{2}\right\rceil+1$,
the $C$-neighbourhood of which is clear and locally admissible,
so by the $C$-reconstruction property, there is a unique translation of $G$
and a unique periodic configuration $\omega_c^\square\in\Omega_\F^\square$
that matches $\omega$ on $I_c\cap (G+k)^c$ (for the right translation).
We can do likewise for any other cell.

Now, two neighbouring cells $c,c'\in\Z^d$
share a common $p$-square window which fixes the same choice of translation for $G$.
Hence, $\omega_c^\square$ and $\omega_{c'}^\square$ overlap on this $p$-square,
and $\Omega_\F^\square$ is $p$-periodic, so they are equal.
Thus, all the cells of the infinite connected component $I(b)$ must agree
on the same $\omega^\square$.
The map $\phi:(\omega,b)\mapsto\omega^\square$
is measurable
(for $\epsilon$ small-enough, so that $I$ has density greater than $\frac{1}{2}$).

In particular, $\omega$ and $\omega^\square$ can only differ outside of $I(b)$,
or on the translation of $G$,
so $d_H\left(\omega,\omega^\square\right)\leq \text{density}(I)+\text{density}(G)\leq 
48\left(2\left(C+\left\lceil\frac{p}{2}\right\rceil\right)+1\right)^d\epsilon+\rho$,
and the same bound holds for $d_B\left(\pi_1^*(\lambda),\phi^*(\lambda)\right)$.
At last, we can fill-in the $\square$ symbols of $G$ in an appropriate random way,
in order to send $\phi^*(\lambda)$ into $\M_\F$,
without changing the bound on $d_B$.
\end{proof}
\end{proposition}

In particular, this proposition gives us a linear
$O(\epsilon)$ bound for the stability of any actually periodic tiling
(which will be $0$-almost periodic, with $G=\emptyset$ and $C$-reconstruction for some $C$).

However, it doesn't apply to the Red-Black tiling from the previous section,
for which we can juxtapose side by side
a Red and a Black $n$-macro-tile at any scale
in a locally admissible way, which breaks the desired quasi-periodicity.

\begin{corollary}[Stability] \label{cor:Stab}
Assume there is a \emph{sequence} of triplets
$\left(p_n,\rho_n,C_n\right)$ for which
Proposition~\ref{prop:BesicovitchBound} applies to $\Omega_\F$.
Then, as soon as $\rho_n\to 0$,
we conclude that $\Omega_\F$ is a stable SFT.
\end{corollary}

\begin{lemma}[Meta Multi-Scale-to-Polynomial Bound]
Consider $D_k=\epsilon\alpha^k+\beta^k$ with $k\in \Z$ and $0<\beta<1<\alpha$.
Denote $\theta=\log_\alpha\left(\frac{1}{\beta}\right)=\frac{-\ln(\beta)}{\ln(\alpha)}>0$.
Then, for any choice $K\in\Z$,
the following bound holds as long as $\epsilon\leq \frac{\theta}{\alpha^{K(1+\theta)}}$:
\[
\min_{k\geq K} D_k \leq \max\left( \sqrt{\alpha},\sqrt{\frac{1}{\beta}}\right)\times
\left( \theta^{\frac{1}{1+\theta}}+(1/\theta)^{\frac{1}{1+1/\theta}}\right)\times
\epsilon^{\frac{\theta}{1+\theta}}.
\]

\begin{proof}
We will later on find the optimal bound on the right assuming $k\in\R$.
Then, by replacing $k$ with the nearest integer we will either increase the power of $\alpha$
by $\frac{1}{2}$ or decrease the one of $\beta$ by $\frac{1}{2}$.
Note that this bound works best under the assumption that $\alpha\approx\frac{1}{\beta}$.
If one is much bigger than the other, we may simply decide on which side we always round $k$,
with an added factor $\alpha$ or $\frac{1}{\beta}$
instead.

Now, consider the parameter
$x:=\alpha^k\in \R^{+*}$.
Thus, $k=\log_\alpha(x)$ so:
\[
\beta^k
=\exp\left(\frac{\ln(x)}{\ln(\alpha)}\times \ln(\beta)\right)
=\exp\left(-\theta \ln(x) \right) = x^{-\theta} .
\]
With this rewriting, $D(x):=\epsilon x+x^{-\theta}$ is much easier to minimise.
Indeed, $D$ can be seen as a positive convex function
that goes to $+\infty$ on $0^+$ and $+\infty$,
hence is minimised when $D'(x) := \epsilon-\theta x^{-\theta-1}=0$,
thus at $x=\left(\frac{\theta}{\epsilon}\right)^{\frac{1}{\theta+1}}$.
Using this value of $x$ in $D$ directly gives us the rest of the expected bound.

Now, for the domain of validity, for us to be ably to round $k$ properly,
we simply require $k=\log_\alpha(x) \geq K$,
which translates as $\epsilon\leq \frac{\theta}{\alpha^{K(1+\theta)}}$.
When the bound doesn't hold, when $K$ is greater that the optimal value,
the optimal choice is simply $D_K$.
\end{proof}
\end{lemma}

\begin{corollary}[Polynomial Stability] \label{cor:PolyStab}
Assume there is a \emph{sequence} of triplets
$\left(p_n,\rho_n,C_n\right)$ for which
Proposition~\ref{prop:BesicovitchBound} applies to $\Omega_\F$.
If $C_n+p_n=O\left(\alpha^{\frac{n}{d}}\right)$ and $\rho_n=O\left(\beta^n\right)$,
then using the previous lemma gives us a polynomial bound
$O\left(\epsilon^{\frac{\theta}{1+\theta}}\right)$ on the speed of convergence, with $\theta=\frac{-\ln(\beta)}{\ln(\alpha)}$.
\end{corollary}

\begin{remark}[]
To illustrate how this framework applies,
let us use it to obtain the polynomial stability
for the enhanced Robinson tiling.

Unlike the usual Robinson tiling,
the enhanced variant enforces alignment of neighbouring macro-tiles in a local way.
At the scale of $N$-macro-tiles,
if we forget about the grid around these tiles,
of density $\rho_n=1-\frac{\left(2^{n}-1\right)^2}{4^n}=O\left(\frac{1}{2^n}\right)$,
we obtain a $p_n$-periodic tiling with $p_n=2\times 2^n$.
What's more, we can prove the tiling has $C_n$-reconstruction~\cite[Proposition 7.7]{GaySa21},
with a radius $C_n=2^n-1$.
As we have $C_n+p_n=O\left(2^n\right)=O\left(4^\frac{n}{2}\right)$,
we can apply the previous corollary with parameters $(\alpha,\beta)=\left(4,\frac{1}{2}\right)$,
so $\theta=\frac{1}{2}$ and $\frac{\theta}{1+\theta}=\frac{1}{3}$.
Hence, we fall back on the $O\left(\sqrt[3]{\epsilon}\right)$ bound
of the previous article~\cite[Theorem 7.9]{GaySa21} (with a comparable multiplicative constant)
which is to be expected as we basically generalised the scheme of the proof
used in that paper.
\end{remark}

More generally, in a hierarchical tiling,
at the scale of ``macro-tiles'' of diameter $x$,
the typical reconstruction radius we may hope for is of order $x$ at least
(\emph{i.e.}\ the size of a macro-tile), and likewise for the quasi-periodicity.
Conversely, among the $x^d$ cells in a macro-tile, we may have to ignore at least 
a one-dimensional ``wire'' that crosses the whole macro-tile, hence hence $\rho$
of order $\frac{1}{x^{d-1}}$ at least.
Following the same general computations as in the previous lemma,
we conclude that in $d$ dimensions, the best speed of convergence we may obtain is
$O\left(\epsilon^{\frac{d-1}{2d-1}}\right)$.
With $d=2$, we have $\frac{d-1}{2d-1}=\frac{1}{3}$, the order of convergence obtained
for the enhanced Robinson tiling.
The question of whether we can obtain a faster bound for the convergence speed of aperiodic tilings,
whether by improving upon the minimal values of $(C,p,\rho)$ conjectured here
(and in particular on the $C$-reconstruction),
or by using another method altogether, is still open.


\section{Undecidability of the Stability} \label{sec:UndecidableStability}

In the previous sections, we showed how a simple bi-coloured tiling can be unstable,
and how a class of well-behaved quasi-periodic SFTs can be stable.
We will now make full use of these ideas in order to equate the notion of stability
with some undecidable problems in the arithmetical hierarchy
through the emergence of said unstable structure.

Matter-of-factly, proving $\Pi_2$-hardness would directly imply
the weaker bounds we introduce first.
However, the $\Pi_2$-hard construction relies on the $\Sigma_1$-hard one, 
and we believe the $\Pi_1$-hard one uses a complementary and more intuitive idea
that will help get the point across.

\subsection{\texorpdfstring{$\Pi_1$}{Π₁}-hard construction}
\label{subsec:P1hard}

First, we will make use of the halting problem $P_{halt}$.
What we want to do here is to encode computations into the Robinson tiling
in a \emph{stable} way, and make an \emph{unstable} phase emerge
\emph{iff} the machine terminates.
This will equate the $\Sigma_1$-complete halting problem with instability among a class of SFTs,
hence $\Pi_1$-hardness of $P_{stab}$ in general.

In the previous Red-Black example of Section~\ref{sec:RedBlackRobinson},
the main ingredient allowing instability
was the existence of two kinds of $n$-macro-tiles at any scale
(widely different for the finite Hamming distance) instead of just one
(four similar tiles, up to the orientation of their low-density central cross) in the monochromatic case.
The two kinds of macro-tiles cannot coexist in the same generic Robinson configuration,
but we can replace one with the other for a small price in the presence of noise.

\subsubsection{Description of the Tileset}
\label{subsec:TuringRobinson}

Let us first describe the tileset used in this section.
We won't explain in details how Turing machines can be implemented
inside the Robinson tiling, but the interested reader may look at the original article
by Robinson~\cite{Rob71} or lecture notes by Jeandel and Vanier~\cite{JeanVa20}
for a formal study of this simulation result.

We will use a variant of this construction more suited to our needs, with two layers,
defined on the alphabet $\A\subset\A_R\times \A_M$,
where $\A_R$ stands for the common Robinson layer,
and $\A_M$ for the layer specific to a given Turing machine $M$.
Consequently, we will denote $\Omega_{P_1(M)}$ the corresponding SFT.

\begin{figure}[H]
\figaligntop
\centering
\includegraphics[width=.6\textwidth]{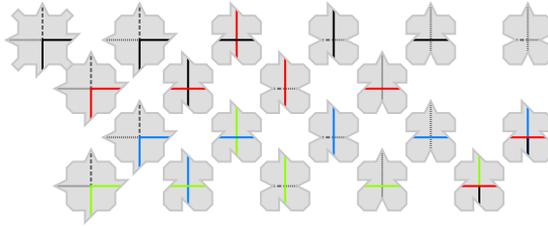}
\figalignmid
\caption{Main tiles of the alphabet $\A_R$.} \label{fig:RobinsonForTuring}
\figalignbot
\end{figure}

Let's first describe the common layer $\A_R$.
As we can see in Figure~\ref{fig:RobinsonForTuring},
the tileset uses four main colours, as well as grey dotted and dashed lines.
These grey lines must match with one of the same type (either dotted \emph{or} dashed on both
sides of an edge),
and serve to enforce alignment of the Robinson 
macro-tiles locally,
to guarantee stability of the structure itself,
just like for the enhanced Robinon tiling~\cite[Proposition 7.7]{GaySa21}.
Notice how bumpy-corners must be Black,
after which we alternate between Black and Red.
At some point, to-be-decided by the layer $\A_M$,
we may transition from the Red-Black (stable) regime to the Blue-Green (unstable) regime
using one of the two transition tiles on the bottom-right of Figure~\ref{fig:RobinsonForTuring}.
The whole set $\A_R$ is given by all the rotations of the first three columns
(but no symmetry, so that we may preserve the chirality of macro-tiles,
so that each arm of the central cross may indicate 
the overall orientation of the macro-tile)
and rotations and symmetries of the rest,
which brings us to a total of $\left|\A_R\right|=172$ tiles.

\begin{figure}[H]
\figaligntop
\centering
\includegraphics[width=.5\textwidth]{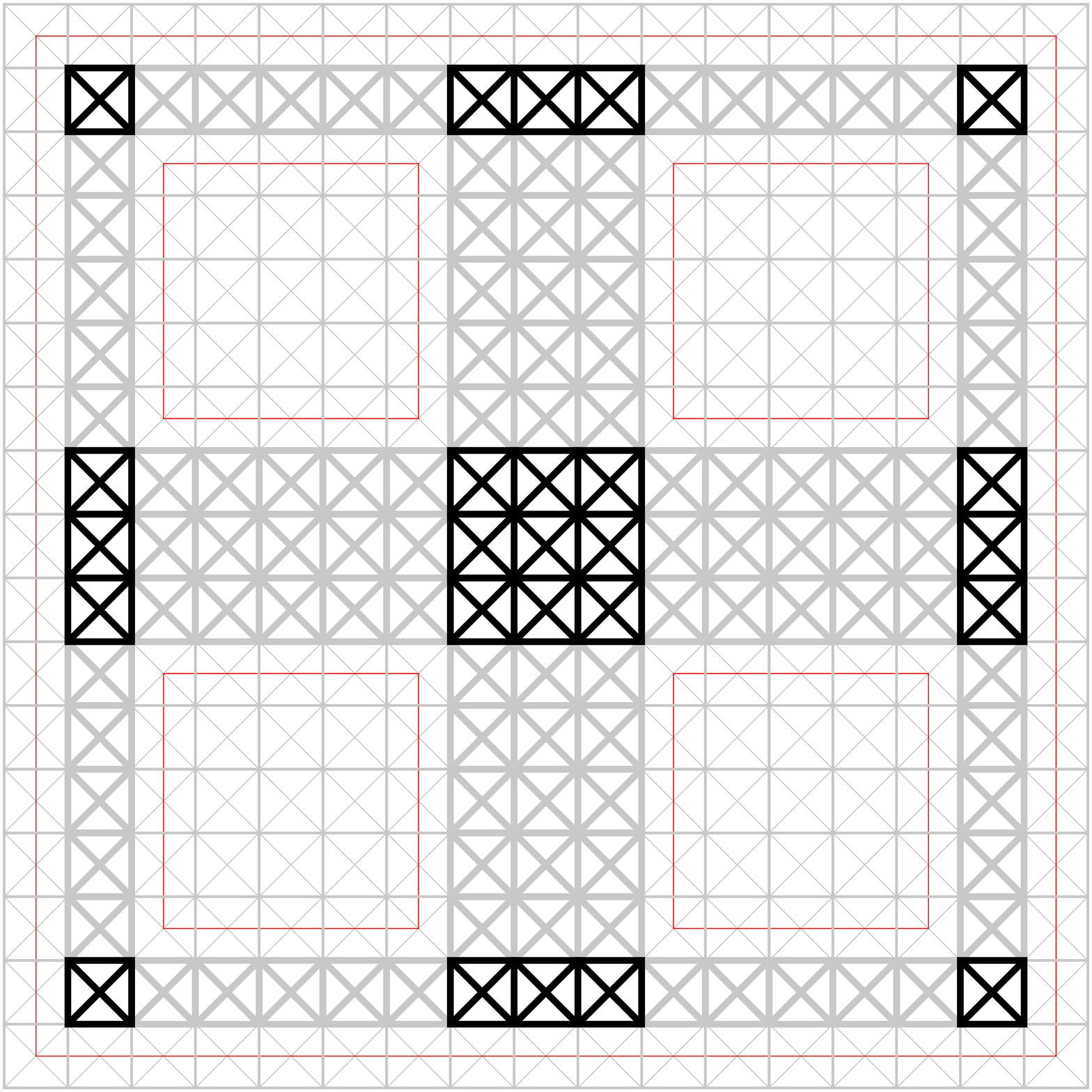}
\figalignmid
\caption{Space-time diagram of a Turing machine\\in a $5$-macro-tile.} 
\label{fig:SpaceTimeDiagram}
\figalignbot
\end{figure}

Now, without detailing the intricacies of $\A_M$ and how it is coupled with $\A_R$ in $\A$,
let us give the general idea and specificities of our construction.
Here, each Red square (of length $4^n+1$, in the center of a $(2n+1)$-macro-tile)
will contain a limited space-time diagram of the Turing machine $M$ with a semi-infinite ribbon,
while avoiding smaller red squares which contain their own space-time diagram.
This is illustrated in Figure~\ref{fig:SpaceTimeDiagram},
where the black crossed cells represent
the patches of space-time diagram,
and the grey cells are communication channels
that synchronise the otherwise disconnected
patches of the diagram.
The $n$-th scale of simulation, occurring in a $(2n+1)$-macro-tile,
thus has a space-time horizon of $2^n+1$ tiles,
initiated on the empty input on the bottom row.

The main difference with the canonical construction is how it behaves when $M$ stops.
In Robinson's article, the tiling doesn't allow for $M$ to stop,
in order to prove that the tileability problem is undecidable.
Here, when $M$ halts in the $n$-th scale of simulation,
it idles until the border of the square can ``notice'' the halting,
and decide freely whether it will force a transition from Black to Blue or Green on its border.
After which, at higher scales,
no more computations occur.
Still, whether or not this transition occurs, we have arbitrarily big macro-tiles,
and thus $\Omega_{P_1(M)}\neq\emptyset$.
Note that a description of $M$ can be algorithmically converted into
the set of forbidden patterns $P_1(M)$
(and its corresponding alphabet) in finite time.

\begin{theorem}[$P_{stab}$ is $\Pi_1$-hard]
Consider a Turing machine $M$.
Then the SFT defined by $P_1(M)$ is stable ($P_1(M)\in P_{stab}$)
\emph{iff} $M$ does \emph{not} halt on the empty input ($M\notin P_{halt}$).
As $P_{halt}$ is $\Sigma_1$-complete, we deduce that $P_{stab}$ is $\Pi_1$-hard.
\end{theorem}

The following subsubsections will each focus on one of the implications,
which put together directly give the previous result.

\subsubsection{The Stable Case}

For the stable case, assume that $M\notin P_{halt}$.
Because of this, at any scale of admissible macro-tiles,
the previously described transition from the Red-Black to the Blue-Green regime cannot occur,
and the two last lines of tiles in Figure~\ref{fig:RobinsonForTuring} may as well not exist in $\A_R$.
Our goal is to prove that the framework of Section~\ref{sec:AperiodicStability} applies here.

Notice how we can project the alphabet $\A$
onto its first coordinate $\A_R$
and then erase the information on
which of the four colours is used
for the lines atop of a tile.
This way, we fall back on the enhanced Robinson tiling studied in our previous paper.
In particular, the following structural result 
applies:

\begin{proposition}[{\cite[Proposition 7.7]{GaySa21}}]
Consider the enhanced Robinson SFT.
Let us denote $B_k=\llbracket -k,k\rrbracket^2$.
For any scale of macro-tiles $n\geq 2$, the constant $R_n=2^n-1$ is such that,
for any $k\geq 0$ and any clear locally admissible pattern $\omega$ on $B_{k+R_n}$,
its restriction $\omega|_{B_k}$ is made of
well-aligned and orientated $n$-macro-tiles,
plus the grid around them which we do not control.
\end{proposition}

Thence,
at the $n$-th scale of simulation (\textit{i.e.}\ in $(2n+1)$-macro-tiles),
the tiling $\Omega_{P_1(M)}$ is $\frac{1}{4^n}$-almost $p_n$-periodic (with $p_n=4^{n+1}$)
with $C_n$-reconstruction ($C_n=R_{2n+1}$) \emph{if we specifically look} at the layer $\A_R$.
However, we need to tread a bit more carefully to obtain the desired periodic behaviour on
the other coordinate of the alphabet $\A$,
and we will actually specifically extend the grid around $(2n+1)$-macro-tiles into a larger set $G_n$ to do so.

\begin{lemma} \label{lem:AlmostPeriodicTuringRobinson}
Using the constant choices from the previous paragraph,
the SFT $\Omega_{P_1(M)}$ is $\rho_n$-almost $p_n$-periodic with $C_n$-reconstruction,
with $G_n$ the area outside of Red squares
up to the $n$-th scale and $\rho_n$ its density.

\begin{proof}
Because $\Omega_{P_1(M)}$ enforces alignment
in a local way,
for \emph{any} tiling $\omega$
we obtain the same set $G_n$ (up to translation)
by looking at all the tiles outside of Red squares
up to the $n$-th scale of simulations.
This $G_n$ is $p_n$-periodic,
and in particular includes
the grid surrounding $(2n+1)$-macro-tiles so that
we have $\rho_n$-quasi $p_n$-periodicity
with $C_n$-reconstruction
\emph{on the layer} $\A_R$.

Regarding alignment, notice that
$G_n$ has the same periodicity as
the grid around $(2n+1)$-macro-tiles,
whose alignment
is fixed by the $C_n$-reconstruction
on the layer $\A_R$,
hence its alignment is fixed in the same way.

Remark that, on the layer $\A_M$,
because the Turing machine is deterministic,
everything that happens on the \emph{inside} of a given admissible Red square is fixed,
insulated from outside interference.
Hence, on this layer
(and using of course the alignment of $G_n$
given by the layer $\A_R$)
we obtain a $\frac{p_n}{2}$-periodic behaviour outside of $G_n$,
as it does not depend on the orientation of the $(2n+1)$-macro-tiles.
\end{proof}
\end{lemma}

In order to conclude, we need to compute $\rho_n$ the density of $G_n$.

\begin{lemma}
In a $(2n+1)$-macro-tile, we have $O\left(12^n\right)$ tiles outside
of the Red squares.

\begin{proof}
The general idea of the proof is that Red squares
form a kind of Sierpiński carpet inside macro-tiles.

Denote $r_n$ the number of tiles \emph{inside} the Red squares in a $(2n+1)$-macro-tile.
As we can see on Figure~\ref{fig:RBStruct},
in the process of forming a $(2n+3)$-macro-tile,
we will create a big central square around four $(2n+1)$-macro-tiles,
surrounded by twelve $(2n+1)$-macro-tiles.
As we already know the size of this big square, we obtain the following recurrence:
\[
r_{n+1}=12 r_n + \left(4^{n+1}+1\right)^2 \geq 12 r_n+16^{n+1} .
\]
As $r_1=25\geq 16$,
we obtain by induction $r_n\geq 4^{n+1}\left(4^n-3^n\right)$.
At the same time, a $(2n+1)$-macro-tile
has $\left(2^{2n+1}-1\right)^2\leq 4^{2n+1}$ tiles in total,
so \emph{at most} $4^{n+1}3^n=4\times 12^n$ tiles outside the Red squares.
\end{proof}
\end{lemma}

Hence, as $(2n+1)$-macro-tiles use $\Theta\left(16^n\right)$ tiles in total,
we conclude that $G_n$ has density $\rho_n=O\left(\left(\frac{3}{4}\right)^n\right)$.

\begin{proposition}[] \label{prop:P1CvRate}
Consider $M\notin P_{halt}$.
Then $\Omega_{P_1(M)}$ is polynomially stable,
with convergence speed $O\left(\epsilon^r\right)$ at rate
$r=\frac{2-\log_2(3)}{6-\log_2(3)}\approx 0.094$.

\begin{proof}
We apply Corollary~\ref{cor:PolyStab},
with constants $\alpha=4$ and $\beta=\frac{3}{4}$,
so $r=\frac{\theta}{1+\theta}$ gives the announced rate.
\end{proof}
\end{proposition}

\subsubsection{The Unstable Case}

\begin{proposition}[] \label{prop:P1colorflip}
Assume $M\in P_{halt}$.
Then for any $\epsilon>0$ we have a measure $\mu\in\M_{P_1(M)}^\B(\epsilon)$
such that $d_B\left(\mu,\M_{P_1(M)}\right)\geq \frac{1}{4^{n+1}}$,
where $n$ denotes the last scale of simulation, at which $M$ halts.

\begin{proof}
Consider $N=2(n+1)$ the first scale at which $N$-macro-tiles
have a big Blue or Green square in the middle.
Assuming two aligned $N$-macro-tiles don't use the same colour for the square
(of diameter $d=2^{N-1}+1$ tiles),
then we obtain at least $p=4\times (d-1)=2^{N+1}$ differences.

By following the very same colour-flipping process
as in Proposition~\ref{prop:UnstableMeasure},
but on the Blue-Green bit starting at the scale of $N$-macro-tiles,
we obtain a generic colour-flipped configuration $\omega$
(with monochromatic Blue or Green squares in the $N$-macro-tiles).

Thus, for any generic $\omega\in\Omega_{P_1(M)}$ that aligns with $\omega'$
up to the scale of $N$-macro-tiles, 
we obtain a lower bound
$d_H\left(\omega,\omega'\right) \geq \frac{1}{2}\times \frac{p}{4^N}=\frac{1}{2^N}=
\frac{1}{4^{n+1}}$,
with the factor $\frac{1}{2}$ coming from the frequency of Blue and Green big squares in
$\omega'$, whereas \emph{all} such squares of $\omega$ must be of the same colour.

Now, assume that $N$-macro-tiles in $\omega$ and $\omega'$ don't align well.
By choosing the best pairing of $N$-macro-tiles between $\omega$ and $\omega'$,
we still have a rectangle with both sides of length at least $2^{N-1}-1$
(the size of a $(N-1)$-macro-tile) where the $N$-macro-tiles of both tilings overlap.
In this area, both macro-tiles have a Blue or Green corner of their big square,
made of at least $2\times\left(2^{N-1}+1\right)$ tiles.
As these two corners intersect in at most $2$ tiles,
and the rest of the area is guaranteed to use only Black or Red communication channels,
we have at least $2^N$ differences between $\omega$ and $\omega'$ in this window.
As this process repeats $2^N$ periodically in both directions,
without even having to take the colour-flipping into account,
we still obtain $d_H\left(\omega,\omega'\right)\geq \frac{1}{2^N}=\frac{1}{4^{n+1}}$.
\end{proof}
\end{proposition}

\begin{remark}[]
More generally, as long as we can guarantee \emph{one}
difference between the two kinds of macro-tiles which we colour-flip,
we obtain a lower bound on $d_B$ of order $\frac{1}{\text{tile area}}$.
We will directly invoke this ``obvious'' lower bound for further unstable cases.

Still, the order of magnitude $\frac{1}{\text{tile diameter}}$ obtained in the previous proposition
is the best one can reasonably hope for in general,
as a signal that transits through a macro-tile will typically only cross a number of tiles
proportional to the diameter, normalised by the tile area.
\end{remark}

\subsection{\texorpdfstring{$\Sigma_1$}{Σ₁}-hard Construction}
\label{subsec:S1hard}

We can ``flip around'' the previous construction,
by adding an \emph{unstable} information atop of the structure simulating the Turing machine,
in such a way that the information gets frozen and becomes \emph{stable} if the machine halts.
We will first describe the construction of $S_1(M)$ out of a machine $M$,
and then state the corresponding indecidability result.

In the previous tileset $P_1(M)$, the Robinson layer $\A_R$ used
one communication channel with four different colours.
Here, for $S_1(M)$, we use \emph{two} communication channels in the lines of the Robinson structure,
each one having two possible values.
First, the Red-Black channel must be initialised as Black in bumpy corners, and then alternate,
in order to have the right structure
to simulate the machine $M$.
Second, the Blue-Green channel can be freely initialised.
However, if $M$ halts at a given scale of simulation,
then the border of the Red square \emph{must} be Blue on the other channel,
which we call a \emph{freeze}.
Note that here, we can keep simulating $M$ at higher scales after it halts for the first time,
as subsequent freezes will just occur at scales of macro-tiles where the Blue-Green channel
would be frozen into Blue anyway.

\begin{proposition}[]
We have $S_1(M)\in P_{stab}$ \emph{iff} $M\in P_{halt}$.
Thus, $P_{stab}$ is $\Sigma_1$-hard.

\begin{proof}
First, assume that $M\notin P_{halt}$.
Then we can freely do a colour-flipping process
starting from any $\mu\in\M_{S_1(M)}$,
just like in Proposition~\ref{prop:UnstableMeasure}.
We can start flipping the Blue-Green channel
at the scale of bumpy corners,
hence instability with a $\frac{1}{8}$ lower bound on $d_B$.

Now, assume $M\in P_{halt}$.
Then, in any tiling $\omega\in\Omega_{S_1(M)}$,
the Blue-Green channel is retroactively frozen all the way down to the Green bumpy-corners.
By using the same grid $G_n$ as in Lemma~\ref{lem:AlmostPeriodicTuringRobinson},
we can likewise ignore everything that happens outside of Red squares, and control everything inside,
hence a $\rho_n$-almost $p_n$-periodic tiling with $p_n=O\left(4^n\right)$
and $\rho_n=O\left(\left(\frac{3}{4}\right)^n\right)$.

Finally, denote $n_{halt}$ the first scale of simulation at which $M$ halts in $S_1(M)$.
If we try to reconstruct things locally at steps lower than $n_{halt}$,
then we will reach a family of well-aligned and well-oriented $(2n+1)$-macro-tiles,
but without any freezing happening in the tiles, hence this Blue-Green channel
that may not behave in a globally admissible way,
all the way down to the high-density set of bumpy-corners.
Still, as long as $n\geq n_{halt}$, the freezing prevents this from happening,
and using the same $C_n=O\left(4^n\right)$ as in 
Lemma~\ref{lem:AlmostPeriodicTuringRobinson},
we conclude that this scale of the tiling has indeed $C_n$-reconstruction.

Still, starting at high-enough scales, for low-enough values of $\epsilon$,
the proof of Proposition~\ref{prop:P1CvRate} applies \emph{verbatim},
so we have stability with a polynomial $O\left(\epsilon^r\right)$ convergence rate.
\end{proof}
\end{proposition}

\subsection{\texorpdfstring{$\Pi_2$}{Π₂}-hard Construction}

In the construction for $S_1(M)$, we obtained stability \emph{iff}
there exists a time step such that $M$ halts \emph{on the empty input}.
Consequently, if we manage to twist the construction to include \emph{any possible input},
then we may equate stability with the $\Pi_2$-complete totality problem $P_{total}$.

There are several ways to proceed, but we choose here to use the method
of Toeplitz encoding of the input, because it is quite versatile,
and may more generally be able to convert
a (structurally close to) uniquely ergodic SFT encoding a $\Sigma_k$-hard problem
into a (definitely not uniquely ergodic anymore) SFT encoding a $\Pi_{k+1}$-hard problem.

\subsubsection{Toeplitz Input}

The Toeplitz encoding of an infinite sequence $u\in \Gamma^{\N^*}$ on an alphabet $\Gamma$
consists of inductively filling with $u_n$ half of the holes still free after the previous iterations,
which gives a sequence $u_1*u_1*u_1*u_1*\dots$,
then $u_1u_2u_1*u_1u_2u_1*\dots$, and so on.
Toeplitz sequences have been studied as dynamical systems for a long while~\cite{JaKea69},
and have since been encoded in higher-dimensional SFTs~\cite{CaVa21}.

The idea of the method is to sequentially write the wanted input $u$
into the consecutive scales hierarchical structure,
which will appear as a Toeplitz encoding $u_1u_2u_1\dots$
from the point of view of the simulated Turing machine,
and then adapt the machine to decode it back into its original form $u$ at first.
This method was already used by Barbieri and Sablik~\cite{BarSa19} in particular.

More precisely, we build the tileset $P_2(M)$ as follows.
For the Robinson structure,
we use the same parallel Red-Black and Blue-Green bits as for $S_1(M)$.
We add another channel that can take values in $\Sigma\sqcup\{\#,\$^\Sigma,\$^\#\}$
where $\Sigma$ is the input alphabet of the machine $M$, $\#$ the blank tape symbol,
and the $\$^*$ symbols two supplementary letters.
On Black channels, we can freely
use any symbol $\$^*$
following a letter from $\Sigma$ on the previous Red scale,
but we \emph{must} use $\$^\#$ following $\#$.
On Red channels,
we \emph{must} use a letter from $\Sigma$ following a $\$^\Sigma$ symbol on the previous Black scale,
and use $\#$ following $\$^\#$.
If we look only at the Red channels,
this gives an infinite word $u\in\Sigma^*\#^\N\sqcup\Sigma^\N$.
When $u\in\Sigma^*\#^\N$, we will identify it with its prefix in $\Sigma^*$,
followed by $\#^\N$.

Quite importantly, the choice of a letter is not only communicated
along the regular Red-Black channels in two directions from the center of a macro-tile arm,
but also along the alignment channels of the enhanced Robinson self-aligning structure,
the dotted and dashed lines
in the other two directions.
Thus, any two (well-aligned) neighbouring $N$-macro-tiles must encode the same sequence.

On the simulation layer, the Turing machine is able to read which symbol is written down in the column on the right of its current position.
Hence, from the point of view of the Turing machine simulated in a Red square,
this represents a \emph{read-only} second tape.
In order to adequately use $u$ as an input, we first need to explain what the machine sees.

\begin{lemma}[Toeplitz Encoding of the Input]
Let $u\in\Sigma^*\#^\N\sqcup\Sigma^\N$.
Define $w_n=w_{n-1}u_nw_{n-1}$ by induction,
initialised with the empty word $w_0$.
The word $w_n$ is a prefix of the Toeplitz encoding of the whole sequence $u$.

At the $n$-th scale of simulation, from the point of view of the Turing machine,
the read-only tape reads as $w_{n-1}\$^*\$^*w_{n-1}u_n$.

\begin{proof}
The last letter of the read-only ribbon obviously correspond to the right border of the $n$-th Red square,
hence reads as $u_n$.
The central $\$^*$ symbols come
from the fact, as highlighted by the blue columns in Figure~\ref{fig:Toeplitz},
they correspond to the $(n+1)$-th scale for Black squares
followed by the first scale of bumpy corners.

\begin{figure}[H]
\figaligntop
\centering
\includegraphics[width=.6\textwidth]{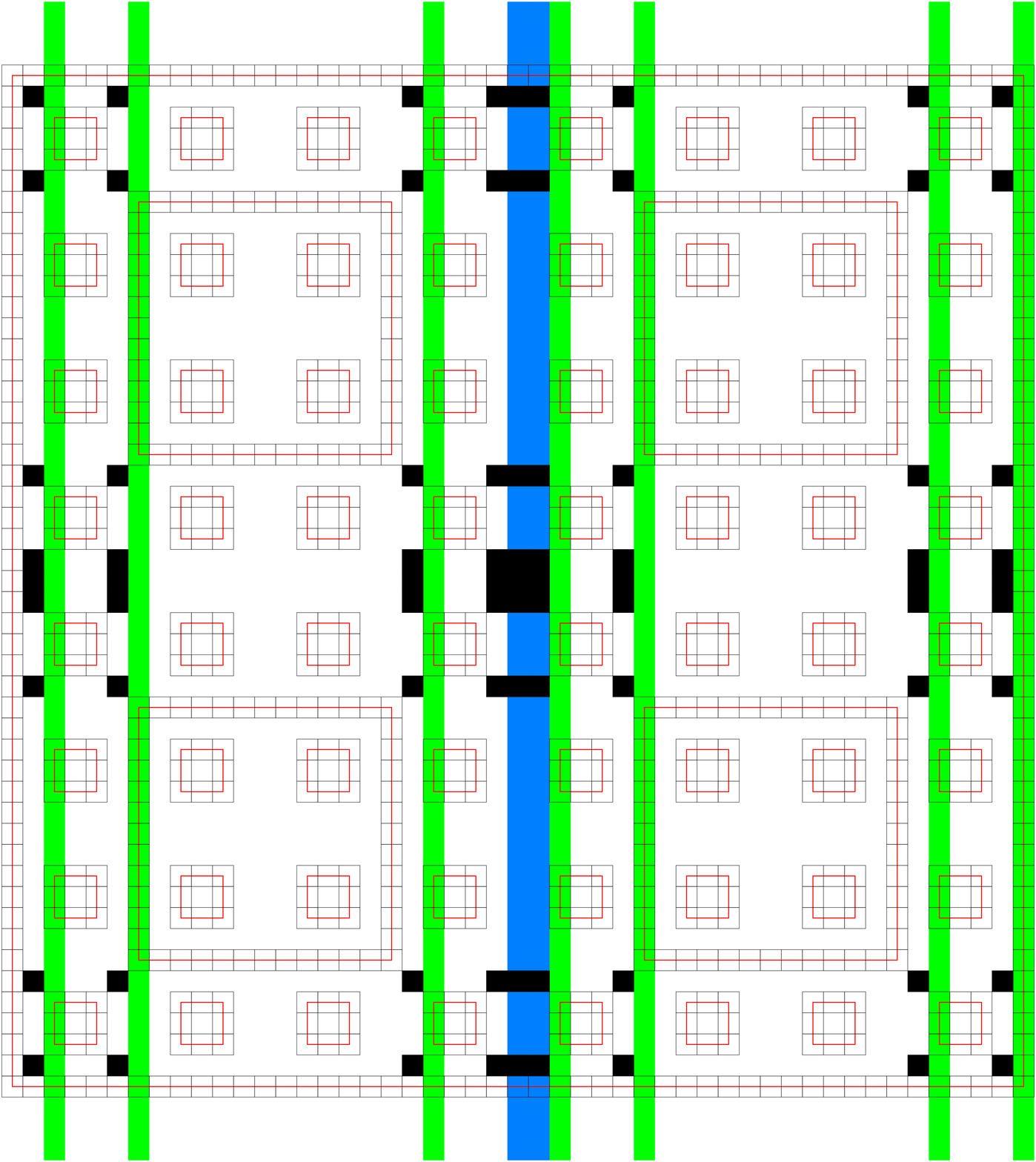}
\figalignmid
\caption{Structure of the read-only input.\\
The highlighted columns are where the read-only values are stored,\\
whereas the machine operates within the black patches.}
\label{fig:Toeplitz}
\figalignbot
\end{figure}

The rest of the word, that reads as $w_{n-1}$
on both ends of this central line,
can be explained by the inductive construction of macro-tiles.
Indeed, each quarter of the $n$-th Red square
is actually a whole $(2n-1)$-macro-tile
with a central Red square,
and the Red squares are themselves stacked in a Toeplitz way
within the macro-tile,
with a gap in-between each that allows to read the letter on them.
\end{proof}
\end{lemma}

\subsubsection{From Decoding the Input to Computations}

Let us explain what Turing machine is encoded into $P_2(M)$,
and how it affects the Blue-Green channel.

First, the machine will have to \emph{decode} the Toeplitz input,
while keeping the Blue-Green channel \emph{stable}
(by using a third non-alternating colour).
More precisely,
the machine will step by step read the letters at positions $2^k$ on the read-only tape
and write them one after another at the beginning of its working tape.
This process will decode the Toeplitz encoding $w_n$ back into the sequence $u_1\dots u_n$.
Using a unary counter, which we multiply by two after reading each letter,
reaching the $k$-th letter will require about $\Theta\left(4^k\right)$ steps of computation.

Now, this process can halt in two ways.
First, we read a $\$^*$ symbol, meaning that we reached halfway through the read-only ribbon.
In this case, the machine simply idles for the rest of its finite runtime,
without unfreezing the Blue-Green bit
when it reaches the top border of the Red square.
This won't happen at big-enough scales of simulation,
considering it would take about $\Theta\left(4^n\right)$ steps but the $n$-th machine only has
a finite horizon of $2^n$ steps,
but it can occur at the initial scales of simulation
and in particular at the very first one where the first symbol is $\$^*$.
Second, we read a $\#$ symbol before reaching the central $\$^*$, in which case the decoding of the word $u\in\Sigma^*$ is complete.
Without waiting, the machine starts then simulating $M$ on $u$
(this will occur roughly at the $2|u|$-th scale of simulation).
This will signal the Red square to \emph{ignite} the unstable Blue-Green bit
(if it was not already done at a lower scale),
as was the case for $P_1(M)$ in Subsection~\ref{subsec:P1hard}.
Then, if $M$ halts on $u$, this will signal the Red square to \emph{freeze} the Blue-Green bit,
as was the case for $S_1(M)$ in Subsection~\ref{subsec:S1hard}.

\subsubsection{Undecidability of the Stability}

\begin{lemma}[]
Assume $M\notin P_{total}$ does not halt on the input $u\in\Sigma^*$.
Consider $\mu_u\in\M_\F$ an invariant measure
with $u\#^\N$ written in the Red scales of any generic configuration.
Then, by colour-flipping the Blue-Green channel after the inital decoding scales, we obtain the measure
$\mu_u^\epsilon\in\M_\F^\B(\epsilon)$,
such that $\inf\limits_{\epsilon>0} d_B\left(\mu_u^\epsilon,\M_\F\right)>0$.

\begin{proof}
As in Proposition~\ref{prop:P1colorflip},
if we compare two macro-tiles with a Blue or Green square,
corresponding to the same input $u\#^\N$,
we can obtain a lower bound on their density of mismatching purely through the Blue-Green square.

Likewise, if we compare such a macro-tile with a macro-tile corresponding to \emph{another} input,
then they must differ in one of the first $|u|+1$ letters in a Red channel.
In this case, we can also obtain a lower bound independent of $\epsilon$,
even if they are perfectly aligned, using this mismatching letter in the input.
\end{proof}
\end{lemma}

\begin{proposition}[]
Denote $\phi(n)$ the first \emph{scale of simulation} at which,
for any input $u\in\Sigma^{\leq n}$, both \emph{decoding} and \emph{computation} are over.
Assume $M\in P_{total}$, so that $\phi(n)<\infty$.
Then, using the notations of Lemma~\ref{lem:AlmostPeriodicTuringRobinson},
the SFT $\Omega_{P_2(M)}$ is $\rho_n$-almost $p_n$-periodic with $C_{\phi(n)}$-reconstruction.

\begin{proof}
If we follow the same scheme of proof as in the lemma,
then we almost obtain a $\rho_n$-almost $p_n$-periodic SFT with $C_n$-reconstruction
up to \emph{one} detail.

Here, the Red-Black channel and all the computations in Red squares behave deterministically,
so they are fixed for a given input (which synchronises between neighbouring tiles),
but the Blue-Green channel is not.
However, if we exploit the $C_{\phi(n)}$-reconstruction of the Robinson structure, then either:
\begin{itemize}
\item a given Red square isn't done decoding its input, so the Blue-Green bit is still frozen, 
uniquely determined,
\item a given Red square has decoded its input $u$, at a scale of simulation lower than $n$,
which implies $|u|<n$, but the Red square actually fits into a bigger $(2\phi(n)+1)$-macro-tile,
which will terminate its simulation of $M$ on $u$,
thus freeze the Blue-Green bit of this Red square.
\end{itemize}
In both cases, we indeed guaranteed that the area inside Red squares is globally admissible,
hence the $n$-th scale of simulation admits $C_{\phi(n)}$-reconstruction.
\end{proof}
\end{proposition}

In particular, $\rho_n\to 0$ so $\Omega_{P_2(M)}$ is stable according to Corollary~\ref{cor:Stab}.
However, because $\phi$ can be roughly as big as any computable function,
we can't possibly exhibit a good bound to apply Corollary~\ref{cor:PolyStab},
and will not obtain polynomial stability this time.
The next theorem directly follows:

\begin{theorem}[]
Consider a Turing machine $M$.
We have $M\in P_{total}$ \emph{iff} $P_2(M)\in P_{stab}$.
As $P_{total}$ is $\Pi_2$-complete, we deduce $P_{stab}$ is $\Pi_2$-hard.
\end{theorem}

Note how this process doesn't adapt
to translate the $\Pi_1$-hard construction into a $\Sigma_2$-hard one.
In order to do this, we would need an added universal quantifier,
which cannot work if we encode only \emph{one} input at a time in a ground configuration.
Hence, in this case we would need to enumerate the inputs inside the tiling in any case.

\begin{remark}[Alternate Construction for $\Pi_2$-Hardness]
Let us conclude this section by briefly describing another construction
relating $P_{stab}$ to $P_{total}$, this time without having to encode any input.

The main idea is here to stack ignition-freezing blocks onto each other.
In the tiling $P_2'(M)$, we enumerate the words of $\Sigma^*$,
\emph{e.g.}\ following a lexicographical order biased by increasing lengths.
After enumerating a new word $u$ we simulate $M$ on it.
Once this simulation ends, we both freeze the lower scales of the Blue-Green bit and
ignite an independent Blue-Green bit for higher scales.
Then, we enumerate the next word, rinse and repeat.

If $M\notin P_{total}$, then $M$ will never end computing on $u$,
never freeze this Blue-Green block which we will be able to colour-flip.
If $M\in P_{total}$, then at any given scale of simulation,
there exists a higher scale of simulation at which $M$ terminates on \emph{some} word,
which will guarantee the current scale is frozen.
\end{remark}


\section{\texorpdfstring{$\Pi_4$}{Π₄} Upper Bound on the Stability}
\label{sec:UpperBound}

As announced, we will now need to dig deeper into the framework of
computable analysis on measures to describe how much computation power
is actually needed to decide our notion of stability.

The general idea of computable analysis is to study problems relating to real numbers,
involving continuous functions or differential equations for example,
from the point of view of effective computations~\cite{Wei00}.
Here, we are interested in doing computable analysis specifically on \emph{probability measures}.
The topic has already been studied~\cite{GaHoyRo11,HeSa18,Wei17} but,
given the lack of a widespread theory,
we will introduce all the needed notations and keep things self-contained.

For the rest of this section, we will consider an alphabet $\A$ and a set of forbidden patterns $\F$,
without any more assumptions such as $\Omega_\F\neq \emptyset$.
From there, our goal in this section is to explain a process
to conclude on whether $\F\in P_{stab}$ or not.

Even though we are interested in convergence for the Besicovitch distance $d_B$,
we will actually need to use the weak-* topology.
Indeed, this topology admits an explicit countable basis dense in the set of
$\sigma$-invariant measures on $\A$ (\emph{i.e.}\ the full-shift $\M_\A$),
which is the bedrock upon which most of computable analysis relies.
Hence, before doing anything meaningful with this topology,
we will first introduce our notations to work with it
in Subsection~\ref{subsec:WeakTopo},
and in particular the (family of)
computable distances we will use later on.

Once this preliminary work is done,
we will see how the measure sets $\M_\F$,
$\M_\F^\B(\epsilon)$ and $J\left(\mu,\M_\F\right)$ can be described in this framework.
At last, we will use these descriptions to prove
a $\Pi_4$ upper bound on the problem $P_{stab}$.

\subsection{A Crash Course in Computable Analysis} \label{subsec:WeakTopo}

\begin{definition}[Weak-* Topology]
The weak-* topology on a set of probability measures $\Prob\left(\A^X\right)$
with $\A$ finite and $X$ countable is defined as follows.
We have the convergence $\mu_n\overset{*}{\longrightarrow}\mu$ when,
for any finite subset $U\Subset X$, and any pattern $w\in\A^U$,
we have the convergence $\mu_n([w])\to\mu([w])$.
\end{definition}

\begin{definition}
The weak-* topology is metrisable, induced notably by:
\[
d_r^+(\mu,\nu):=\sum\limits_{n\in\N}\frac{1}{2^n}\times 
\frac{1}{r^{\left|U_n\right|}}\times\frac{1}{ \left|\A^{U_n}\right|}\sum\limits_{w\in\A^{U_n}}
\left|\mu([w])-\nu([w])\right| ,
\]
with $\left(U_n\right)_{n\in\N}$ an increasing sequence of sets that covers the space $X$,
and $r\geq 1$ a normalisation factor (with the convention $d^+=d_1^+$).
In particular, when $X=\Z^d$, we can take $U_n=B_n$.
$\Prob\left(\A^X\right)$ is a compact space for this topology.
When $X=\Z^d$, the space of $\sigma$-invariant measures $\M_\A$ is a closed subspace.
In this subspace, we can and will instead use $U_n:=\llbracket 0,n \rrbracket^d$ to define $d^+$ for the rest of this article.
\end{definition}

Note how, if $V_1\subset \M_{\A_1}$ and $V_2\subset\M_{\A_2}$ are both weakly closed,
then so is the set of their joinings $J\left(V_1,V_2\right)\subset\M_{\A_1\times\A_2}$.

\begin{definition}[Closed Ball]
We denote the closed ball around $\mu\in \M_\A$ of radius $\epsilon >0$ as
$\Ball(\mu,\epsilon)=\left\{\nu\in\M_\A, d^+(\mu,\nu)\leq\epsilon \right\}$,
and this definition extends to the $\epsilon$-neighbourhood of any set of measures.
\end{definition}

\begin{definition}[Periodic Measure]
For $w\in \A^{U_n}$, we denote $w^{\Z^d}\in\Omega_\A$
the configuration obtained through periodic repetition of $w$ in each direction,
and then $\meas{w}= \frac{1}{\left|U_n\right|}\sum\limits_{k\in U_n}
\delta_{\sigma_k\left(w^{\Z^d}\right)}$ the corresponding $\sigma$-invariant measure.
We call such measures \emph{periodic}.
\end{definition}

\begin{lemma}[Covering Lemma for $\M_\A$] \label{lem:CoveringMA}
There is $\psi:\Q^{+*}\times\N^2\to\N$ a computable map such that,
for any finite alphabet $\A$, any dimension $d$ and any rational $\delta>0$:
\[
\M_\A = \bigcup_{w\in\A^{U_{\psi(\delta,|\A|,d)}}}\Ball\left(\meas{w},\delta\right) .
\]

\begin{proof}
Denote $s_n$ the partial sum up to rank $n$
associated to the distance $d^+$.
We can bound $r_n=d^+-s_n\leq \frac{1}{2^{n-1}}$
independently of the dimension $d$, of $\A$ and of any pair of measures.
Hence, for a given value of $\delta$,
we first compute $n(\delta)=2+\left\lceil \log_2\left(\frac{1}{\delta}\right)\right\rceil$,
such that $r_n\leq \frac{\delta}{2}$.
We now want to cover $\M_\A$
with balls of radius $\frac{\delta}{2}$
for the pseudo-distance $s_n$.

Notice that for any $k\leq n$ and any word $w\in\A^{U_k}$,
we have the decomposition
$\mu([w])=\sum_{v\in\A^{U_n},\,v|_{U_k}=w}\mu([v])$.
It follows that:
\[
\begin{array}{rcl}
s_n(\mu,\nu)&=&
\sum\limits_{k\leq n} \frac{1}{2^k\left|\A^{U_k}\right|}
\sum\limits_{w\in \A^{U_k}} |\mu([w])-\nu([w])| \\
&\leq&\sum\limits_{k\leq n}\frac{1}{2^k}\sum\limits_{w\in \A^{U_k}} |\mu([w])-\nu([w])| \\
&\leq&\sum\limits_{k\leq n}\frac{1}{2^k}\sum\limits_{w\in \A^{U_n}} |\mu([w])-\nu([w])| \\
&\leq& 2\sum\limits_{w\in \A^{U_n}} |\mu([w])-\nu([w])| \\
&\leq& 2\left|\A^{U_n}\right| \sup\limits_{w\in \A^{U_n}} |\mu([w])-\nu([w])|.
\end{array}
\]

Hence, we now need to uniformly approximate any $\mu\in\M_\A$
on the window $U_n$ by a periodic measure to conclude.

To do so, consider $\mu_m$ the restriction of $\mu$ to $U_m$.
We identify $\mu_m$, a measure on $\A^{U_m}$, with the measure on $\Omega_\A$
that charges a periodic word $w^{\Z^d}$ (with $w\in\A^{U_m}$)
with probability $\mu([w])=\mu_m(\{w\})$.
Remark that $\mu_m$ is not $\sigma$-invariant, but is $m\Z^d$-periodic under translation,
so we can define the corresponding averaged measure
$\widehat{\mu_m}:=\sum_{w\in\A^{U_m}} \mu_m(\{w\})\times\meas{w}$
which is $\sigma$-invariant.

In particular for any $w\in\A^{U_n}$,
as long as $k+U_n\subset U_m$ (\textit{i.e.}\ $k\in U_{m-n}$),
then $[w]_k:=\sigma_k([w])$ is still a cylinder defined inside $U_m$.
Hence, for any such translation we have $\mu_m\left([w]_k\right)=\mu\left([w]_k\right)=\mu([w])$.
Now, we have:
\[
\widehat{\mu_m}([w]):=\frac{1}{\left|U_m\right|}\sum_{k\in U_m} \mu_m\left([w]_k\right)
= \frac{\left| U_{m-n}\right|}{\left|{U_m}\right|} \mu([k])
+\frac{1}{\left|{U_m}\right|}\sum\limits_{k\in U_m\backslash U_{m-n}} \mu_m\left([w]_k\right),
\]
hence $\widehat{\mu_m}([w])=\mu([w])+O\left(\frac{n}{m}\right)$,
where the computable domination constant depends on $d$.
Now, if we use instead $\mu_m^k$ a dyadic approximation of $\mu_m$ on $\A^{U_m}$,
with precision $\frac{1}{2^k}$,
we obtain a measure $\widehat{\mu_m^k}$ for which:
\[
\widehat{\mu_m^k}([w])=\mu([w])+O\left(\frac{n}{m}\right)
+O\left(\frac{\left|\A^{U_m}\right|}{2^k}\right).
\]
This new term simply uses the domination constant $1$.
Remark in particular that there is only a finite amount of such dyadic measures
on the window $U_m$ with precision $\frac{1}{2^k}$.
We just need to be able to approximate these by periodic measures to conclude.

We can decompose any such dyadic measure as
$\widehat{\mu_m^k}=\frac{1}{2^k}\sum_{w\in\A^{U_m}} p(w)\meas{w}$
with weights $p(w)\in\N$ that sum to $2^k$.
Consider now $M=(m+1)\times 2^k-1$.
On the corresponding window $U_M$,
we can fit a total of $2^k$ slices,
each made of windows $U_m$ stacked in all directions but one.
In $p(w)$ such consecutive slices, we write $w$ in each box $U_m$.
This gives us a configuration $\overline{w}\in\A^{U_M}$ such that,
for any $w\in\A^{U_n}$,
we have $\meas{\overline{w}}([w])=\widehat{\mu_m^k}([w])+O\left(\frac{n}{m}\right)$,
once again with a computable domination constant that depends on $d$, so that:
\[
\left|\meas{\overline{w}}([w])-\mu([w])\right|=O\left(\frac{n}{m}\right)
+O\left(\frac{\left|\A^{U_m}\right|}{2^k}\right).
\]

Thus, we can actually compute integers $m(\delta,|\A|,d)$ and $k(\delta,|\A|,d)$ such that
$\left|\meas{\overline{w}}([w])-\mu([w])\right|\leq\frac{\delta}{2}
\times\frac{1}{2\left|\A^{U_n}\right|}$, which we can replace in the supremum bound for $s_n$.
At last, we proved that there exists a pattern
$w\in\A^{\psi(\delta,|\A|,d)}$ such that $\mu\in\Ball\left(\meas{w},\delta\right)$, 
with $\psi(\delta,|\A|,d):=(m(\delta,|\A|,d)+1)2^{k(\delta,|\A|,d)}-1$ a map that can be computed by a Turing machine.
\end{proof}
\end{lemma}

In particular, the density of the family 
$\left( \meas{w},w\in\A^{U_n},n\in\N\right)$ of all periodic measures
directly follows from the Covering lemma.
Note how we always have $d_r^+\leq d_1^+$ when $r\geq 1$,
so the previous Covering lemma more generally applies for all these distances.

Let us conclude this subsection with a technical lemma
that relates weak distances when projecting.

\begin{lemma}[Projection Lemma] \label{lem:Projection}
Consider two measures $\lambda,\lambda'\in\M_{\A_1\times\A_2}$.
Then $d_{\left|\A_2\right|}^+\left(\pi_1^*(\lambda),\pi_1^*\left(\lambda'\right)\right)\leq
d^+\left(\lambda,\lambda'\right)$.

\begin{proof}
We have:
\[
\begin{array}{rl}
&d_{\left|A_2\right|}^+\left(\pi_1^*(\lambda),\pi_1^*\left(\lambda'\right)\right)\\
=&\sum\limits_{n\in\N}\frac{1}{2^n}\times\frac{1}{\left|\A_2^{U_n}\right|}\times
\frac{1}{\left|\A_1^{U_n}\right|} \sum\limits_{w_1\in\A_1^{U_n}}
\left|\pi_1^*(\lambda)\left(\left[w_1\right]\right)-
\pi_1^*\left(\lambda'\right)\left(\left[w_1\right]\right)\right| \\
=&\sum\limits_{n\in\N}\frac{1}{2^n}\times\frac{1}{\left|\A_2^{U_n}\right|}\times
\frac{1}{\left|\A_1^{U_n}\right|}\sum\limits_{w_1\in\A_1^{U_n}}\left|
\sum\limits_{w_2\in\A_2^{U_n}}\lambda\left(\left[w_1,w_2\right]\right)
-\lambda'\left(\left[w_1,w_2\right]\right)\right|\\
\leq&\sum\limits_{n\in\N}\frac{1}{2^n}\times\frac{1}{\left|\left(\A_1\times\A_2\right)^{U_n}
\right|}\sum\limits_{\left(w_1,w_2\right)\in\left(\A_1\times \A_2\right)^{U_n}}
\left|\lambda\left(\left[w_1,w_2\right]\right)-\lambda'\left(\left[w_1,w_2\right]\right)\right|\\
= &d^+\left(\lambda,\lambda'\right) ,
\end{array}
\]
\emph{i.e.} the announced bound.
\end{proof}
\end{lemma}

\subsection{Computable Descriptions of Measure Sets}

Now, from the point of view of Turing machines,
the main obstruction to discuss the notion of stability is
that it is not obvious how we should proceed to compute $d_B\left(\mu,\M_\F\right)$.
To do so, we will step-by-step reach a characterisation of the sets
$\M_\F$, $\widetilde{\M_\F^\B}(\epsilon)$ and at last $J\left(\mu,\M_\F\right)$.

\begin{lemma}[Covering Lemma for $\M_\F$] \label{lem:CoveringMF}
Consider $\psi$ the radius-to-scale function obtained in
the Covering Lemma~\ref{lem:CoveringMA} for $\M_\A$.
Assume $\F\subset \A^{U_k}$ for some $k\in\N$. Then:
\[
\M_\F \subset \bigcup\limits_{w\in \W_\F(\rho)} \Ball\left(\meas{w},\rho\right)
\]
where $\W_\F(\rho)\subset\A^{U_{\psi(\rho,|\A|,d)}}$ is the set of
patterns on the window $U_{\psi(\rho,|\A|,d)}$ that contain
at most $\phi(\rho,k,|\A|,d)$ forbidden patterns from $\F$,
with the computable map $\phi(\rho,k,|\A|,d):=\left\lfloor 2^k\times\left|\A^{U_k}\right|\times 
\rho\times\left|U_{\psi(\rho,\dots)}\right|\right\rfloor$.

Note that the set $\W_\F(\rho)$ depends on $k$,$|\A|$ and $d$,
but we hide them from the notation as they are either
directly ``written'' in the computer representation of $\F$
or can be deduced from it.

\begin{proof}
The first covering lemma gives us
$\M_\A=\bigcup_{w\in\A^{U_{\psi(\rho,|\A|,d)}}}\Ball\left(\meas{w},\rho\right)$.
If we prove that $\M_\F\cap\Ball\left(\meas{w},\rho\right)=\emptyset$ whenever
$w$ has more than $\phi(\rho,\dots)$ forbidden patterns, this will conclude the proof.

Consider such a pattern $w\notin\W_\F(\rho)$ and $\mu\in\M_\F$.
For any forbidden pattern $u\in\F\subset \A^{U_k}$,
we have $\mu([u])=0$, thence:
\[
d^+\left(\mu,\meas{w}\right)\geq \frac{1}{2^k}\times\frac{1}{\left|\A^{U_k}\right|}
\sum\limits_{u\in\F} \meas{w}([u]) .
\]
Now, by summing the number of occurrences of each forbidden pattern in $w$,
we obtain at least $\phi(\rho,\dots)+1>
2^k\left|\A^{U_k}\right|\times\rho\times\left|U_{\psi(\rho,\dots)}\right|$.
The last factor is precisely the normalisation constant used to define $\meas{w}$,
so that $d^+\left(\mu,\meas{w}\right)>\rho$, which concludes the proof.
\end{proof}
\end{lemma}

\begin{corollary}
Using the previous notations, we have:
\[
\M_\F=\bigcap\limits_{\rho>0}\bigcup\limits_{w\in\W_\F(\rho)} \Ball\left(\meas{w},\rho\right) .
\]

\begin{proof}
The covering lemma for $\M_\F$ holds for any $\rho>0$,
hence by taking the intersection we directly obtain the inclusion $(\subset)$.

Conversely, consider $\mu\in\bigcap\bigcup(\cdots)$.
There exists a sequence of radii $\rho_n\to 0$ and corresponding patterns $w_n$
such that $d^+\left(\mu,\meas{w_n}\right)\leq\rho_n$.
For any given forbidden pattern $u\in\F$, we have:
\[
\meas{w_n}([u])\leq\frac{\phi\left(\rho_n,\dots\right)
+O\left(\psi\left(\rho_n,\dots\right)^{d-1}\right)}{\left| U_{\psi\left(\rho_n,\dots\right)}\right|}
= O\left(\rho_n\right)+O\left(\frac{1}{\psi\left(\rho_n,\dots\right)}\right)\to 0 ,
\]
with the first bound relating to the occurrences of $u$ within
$w_n$ and the second to the occurrences on the interface
between two square blocks of $w_n^{\Z^d}$.
Thus, for the weak-* limit $\mu$ we have $\mu([u])=0$,
so $\mu\in\M_\F$~\cite[Lemma 2.11]{GaySa21}.
\end{proof}
\end{corollary}

One one hand, the covering lemma tells us that any measure $\M_\F$
can be approximated by some periodic measures with an explicit bound on the speed of convergence.
More precisely, all measures $\mu\in\M_\F$ are $\rho$-close to some $\meas{w}$
some $w\in\W_\F(\rho)$,
but not \emph{all} such $\meas{w}$ are necessarily $\rho$-close to $\M_\F$.
This is roughly correlated to the nuance between
locally admissible and globally admissible tilings.
On the other hand, the corollary tells us that,
while we have no computable bound on the speed of convergence,
we still necessarily converge to measures $\M_\F$,
\emph{i.e.}\ the $\rho$-coverings converge to the set $\M_\F$
in the corresponding Hausdorff topology as $\rho\to 0$.
This will allow us to computationally describe $\M_\F$
as the set of all adherence values of these measures.

We now want to move onto the noisy framework,
to obtain a similar result for $\widetilde{\M_\F^\B}(\epsilon)$,
as a subset of $\M_{\widetilde{\F}}$.
In the following proposition,
we denote $s_n$ the partial sum for $d_{|\A|}^+$ on the space of noises $\Omega_{\{0,1\}}$.
In particular, if we use the computable rank $n(\rho)$ introduced in the proof
of the covering lemma for $\M_\A$,
then we can guarantee $s_{n(\rho)}\leq d_{|\A|}^+\leq s_{n(\rho)}+\rho$.

\begin{lemma}[Covering Lemma for $\widetilde{\M_\F^\B}(\epsilon)$] 
\label{lem:CoveringMFeps}
As in the Covering Lemma~\ref{lem:CoveringMF} for $\M_\F$,
we assume here that $\F\subset \A^{U_k}$.
Note that $\M_{\widetilde{\F}}$ uses the alphabet $\A\times\{0,1\}$,
of cardinality $2|\A|$.
For any $\epsilon\in[0,1]$,
we can refine the covering of $\M_{\widetilde{\F}}$
into a covering of $\widetilde{\M_\F^\B}(\epsilon)$:
\[
\widetilde{\M_\F^\B}(\epsilon) = \bigcap\limits_{\rho>0}
\bigcup\limits_{(w,b)\in\widetilde{\W_\F^\epsilon}(\rho)}\Ball\left(\meas{(w,b)},\rho\right),
\]
where $\widetilde{\W_\F^\epsilon}(\rho)\subset \W_{\widetilde{\F}}(\rho)$ is the subset
for which $s_{n(\rho)}\left(\meas{b},\B(\epsilon)^{\otimes \Z^d}\right)\leq \rho$ holds.

\begin{proof}
As before, we need to prove both inclusions.

Consider first $(\subset)$.
As $\widetilde{\M_\F^\B}(\epsilon)\subset \M_{\widetilde{\F}}$,
we already have the inclusion for any $\rho>0$ \emph{if}
we forget about the new condition on $s_{n(\rho)}$,
using the Covering Lemma~\ref{lem:CoveringMF} for $\M_{\widetilde{\F}}$.
Hence, it suffices to prove that if $(w,b)\in\W_{\widetilde{\F}}(\rho)$
does not satisfy the new condition,
then $d^+\left(\meas{(w,b)},\widetilde{\M_\F^\B} (\epsilon)\right)>\rho$.
This directly follows from the fact that for any $\lambda\in\widetilde{\M_\F^\B}(\epsilon)$,
using the Projection Lemma~\ref{lem:Projection} we have:
\[
d^+\left(\meas{(w,b)},\lambda\right)\geq 
d_{|\A|}^+\left(\meas{b},\B(\epsilon)^{\otimes\Z^d}\right)
\geq s_{n(\rho)}\left(\meas{b},\B(\epsilon)^{\otimes\Z^d}\right)  >\rho .
\]

Conversely, consider $(\supset)$.
For any $\lambda\in\bigcap\bigcup(\cdots) \subset\M_{\widetilde{\F}}$,
we have a sequence $\rho_n\to 0$
and patterns $\left(w_n,b_n\right)\in\widetilde{\W_\F^\epsilon}\left(\rho_n\right)$ 
such that $d^+\left(\meas{\left(w_n,b_n\right)},\lambda\right)\leq \rho_n$.
We just need to prove that $\pi_2^*(\lambda)=\B(\epsilon)^{\otimes\Z^d}$.
This comes from the fact that:
\[
\begin{array}{rcl}
d_{|\A|}^+\left(\pi_2^*(\lambda),\B(\epsilon)^{\otimes\Z^d}\right) &\leq&
d_{|\A|}^+\left(\pi_2^*(\lambda),\meas{b_k}\right) +
d_{|\A|}^+\left(\meas{b_k},\B(\epsilon)^{\otimes\Z^d}\right) \\
&\leq& d^+\left(\lambda,\meas{\left(w_k,b_k\right)}\right) +
d_{|\A|}^+\left(\meas{b_k},\B(\epsilon)^{\otimes\Z^d}\right) \\
&\leq& d^+\left(\lambda,\meas{\left(w_k,b_k\right)}\right) +
s_{n\left(\rho_k\right)}\left(\meas{b_k},\B(\epsilon)^{\otimes\Z^d}\right)+\rho_k \\
&\leq& 3\rho_k .
\end{array}
\]
As $k\to\infty$,
we have $d_{|\A|}^+\left(\pi_2^*(\lambda),\B(\epsilon)^{\otimes\Z^d}\right)=0$,
which concludes the proof.
\end{proof}
\end{lemma}

As long as $\epsilon,\rho\in\Q^{+*}$,
then $(b,\epsilon,\rho)\mapsto s_{n(\rho)}\left(\meas{b},\B(\epsilon)^{\otimes \Z^d}\right)$
is computable. Hence, as for the case of $\M_\F$, this covering lemma tells us \emph{both}
that we \emph{can} mathematically approximate
some $\lambda\in\widetilde{\M_\F^\B}(\epsilon)$ with an explicit bound $\rho$,
and that we have a way of describing $\widetilde{\M_\F^\B}(\epsilon)$
as a set of adherence values of a computable sequence
but without an explicit control on the speed of convergence.

\begin{remark}[Approximating $J\left(\mu,\M_\F\right)$]
Later on, to approximate a joining $\lambda \in J\left(\mu,\M_\F\right)$,
we will consider a periodic measure $\meas{\left(w_1,w_2\right)}$ on $\Omega_{\A\times\A}$
such that $\meas{w_1}$ is close $\meas{w}$,
with $\meas{(w,b)}$ on $\Omega_{\A\times\{0,1\}}$
an approximation of a measure that projects to $\mu\in\M_\F^\B(\epsilon)$
obtained through the Covering Lemma~\ref{lem:CoveringMFeps},
and that $\meas{w_2}$ is close to $\meas{w'}$ obtained through
the Covering Lemma~\ref{lem:CoveringMF} for $\M_\F$.
\end{remark}

\subsection{Equivalent Characterisations of Stability}

Let the measurable event
$\Delta:=\bigcup_{a\neq b\in\A}[(a,b)]\subset \Omega_{\A^2}$.
Because we consider $\sigma$-invariant measures, by an ergodic theorem,
we have $d_B(\mu,\nu)=\inf_{\lambda\in J(\mu,\nu)} \lambda(\Delta)$.
In particular, as $\lambda\mapsto\lambda(\Delta)$ depends only
on a finite window in $\Omega_{\A^2}$,
it is a continuous map for the weak-* topology.

Let us remind what it means for $\F$ to induce a stable SFT:
\[
\forall \delta>0,\exists \epsilon>0, \forall\mu\in\M_\F^\B(\epsilon),d_B\left(\mu,\M_\F\right)
\leq \delta .
\]
Notice how, by monotonicity of the definition, we can restrict this formula by quantifying
$\epsilon$ and $\delta$ over the countable set $\Q^{+*}$ instead.
What's more, using the previous rewriting of $d_B$ through joinings,
we obtain the following characterisation:
\begin{equation} \label{eqn:Stab1}
\forall \delta\in\Q^{+*},\exists \epsilon\in\Q^{+*},
\forall\mu\in\M_\F^\B(\epsilon),
\exists\lambda\in J\left(\mu,\M_\F\right), \lambda(\Delta)\leq\delta .
\end{equation}

Note that, by embedding all the patterns of $\F$ in a big-enough square box $U_k$
and enumerating all the $w\in\A^{U_k}$ that contain at least one forbidden pattern,
We can trivially \emph{compute} $\F'$ such that $\Omega_\F=\Omega_{\F'}$.
As stability with Bernoulli noise is a conjugacy invariant~\cite[Corollary 3.15]{GaySa21},
we can equivalently decide whether $\F'$ is stable instead.
Hence, we will without loss of generality assume that $\F\subset\A^{U_k}$
in the following theorem, so that the covering lemmas may apply.

\begin{proposition}
The SFT $\Omega_\F$ is stable \emph{iff} it satisfies the following formula:
\begin{equation} \label{eqn:Stab2}
\begin{array}{c}
\forall\delta\in\Q^{+*},\exists\epsilon\in\Q^{+*},\forall\rho\in\Q^{+*},\\
\forall \mu\in\M_\F^\B(\epsilon),
\exists \left(w_1,w_2\right)\in\left(\A^2\right)^{U_{\psi\left(\rho,\left|\A^2\right|,d\right)}},\\
\left[ d_{|\A|}^+\left(\meas{w_1},\mu\right)\leq\rho \right] \hspace{-2pt}\wedge\hspace{-2pt}
\left[ d_{|\A|}^+\left(\meas{w_2},\M_\F\right)\leq\rho \right] \hspace{-2pt}\wedge\hspace{-2pt}
\left[ \meas{\left(w_1,w_2\right)}(\Delta)\leq\delta+|\A|^2\rho \right] .
\end{array}
\end{equation}

\begin{proof}
Consider first the implication $(\ref{eqn:Stab1}\Rightarrow\ref{eqn:Stab2})$.
Assuming Equation~\ref{eqn:Stab1} is satisfied,
let us fix $\delta$, $\epsilon$ and $\mu$
such that there exists a joining $\lambda\in J\left(\mu,\M_\F\right)$
for which $\lambda(\Delta) \leq \delta$.
Using the Covering Lemma~\ref{lem:CoveringMA} for $\M_{\A^2}$,
we know that for any rational $\rho\in\Q^{+*}$,
we have some couple $\left(w_1,w_2\right)\in \left(\A^2\right)^{
U_{\psi\left(\rho,\left|\A^2\right|,d\right)}}$
such that $d^+\left(\meas{\left(w_1,w_2\right)},\lambda\right)\leq \rho$.
The first two inequalities in Equation~\ref{eqn:Stab2} follow directly
from the Projection Lemma~\ref{lem:Projection}.
For the third one:
\[
\begin{array}{rcl}
\meas{\left(w_1,w_2\right)}(\Delta) &\leq& \lambda(\Delta)+
\sum_{a,b\in\A}\left|\lambda([(a,b)])-\meas{\left(w_1,w_2\right)}([(a,b)])\right| , \\
&\leq& \lambda(\Delta)+|\A|^2 d^+\left(\lambda,\meas{\left(w_1,w_2\right)}\right) , \\
&\leq& \delta+|\A|^2\rho .
\end{array}
\]
Remark that the consecutive universal blocks $\forall\mu\in\M_\F^\B(\epsilon)$
and $\forall\rho\in\Q^{+*}$ do not depend on each other,
so we can freely reorder them as in Equation~\ref{eqn:Stab2},
which concludes this implication.

Conversely, suppose Equation~\ref{eqn:Stab2} holds
and let us prove $(\ref{eqn:Stab2}\Rightarrow\ref{eqn:Stab1})$.
Likewise, fix $\delta$, $\epsilon$ and $\mu$ for which
the rest of the formula (\emph{i.e.}\ $\forall \rho \dots$) is satisfied.
Consider a sequence $\rho_n\to 0$
and the consequent patterns
$\meas{w_1^n,w_2^n}\in\left(\A^2\right)^{U_{\psi\left(\rho_n,\left|\A^2\right|,d\right)}}$.
Up to extraction of a subsequence,
we can assume that the sequence weakly converges
to a measure $\lambda\in\M_{\A^2}$.
The first inequality of Equation~\ref{eqn:Stab2} gives us $\pi_1^*(\lambda)=\mu$ at the limit.
The second inequality gives us $\pi_2^*(\lambda)\in\M_\F$ (as $\M_\F$ is closed),
so $\lambda\in J\left(\mu,\M_\F\right)$.
At last, the third inequality becomes $\lambda(\Delta)\leq \delta$ at the limit
by continuity of $\lambda\mapsto\lambda(\Delta)$, which concludes the proof.
\end{proof}
\end{proposition}

We now want to use the Covering Lemma~\ref{lem:CoveringMFeps}
for $\widetilde{\M_\F^\B}(\epsilon)$ to replace the $\forall\mu$ block
by a universal block that quantifies over rational numbers instead.

\begin{proposition}
The SFT $\Omega_\F$ is stable \emph{iff} it satisfies the following formula:
\begin{equation} \label{eqn:Stab3}
\begin{array}{c}
\forall\delta\in\Q^{+*},\exists\epsilon\in\Q^{+*},\forall\rho\in\Q^{+*},
\exists\gamma\in\Q^{+*},\gamma\leq \rho,\\
\forall (w,b)\in \widetilde{\W_\F^\epsilon}(\gamma),
\exists w_0\in \W_\F(\rho),\exists \left(w_1,w_2\right)\in\left(\A^2\right)^{
U_{\psi\left(\rho,\left|\A^2\right|,d\right)}},\\
\left[ d_{|\A|}^+\left(\meas{w_1},\meas{w}\right)\leq2\rho \right] \wedge
\left[ d_{|\A|}^+\left(\meas{w_2},\meas{w_0}\right)\leq2\rho \right] \\ \wedge
\left[ \meas{\left(w_1,w_2\right)}(\Delta)\leq\delta+|\A|^2\rho \right] .
\end{array}
\end{equation}

\begin{proof}
Let us prove that $(\ref{eqn:Stab2}\Rightarrow\ref{eqn:Stab3})$.
Assume Equation~\ref{eqn:Stab2} holds true and
fix $\delta$, $\epsilon$ and $\rho$ so that the rest of the formula holds true.
The Covering Lemma~\ref{lem:CoveringMFeps} for $\widetilde{\M_\F^\B}(\epsilon)$
implies that
$\max_{(w,b)\in \widetilde{W_\F^\epsilon}(\gamma)} d^+\left(\meas{(w,b)},
\widetilde{\M_\F^\B}(\epsilon)\right)\underset{\gamma\to 0}{\longrightarrow}0$.
In particular, there is a rational $\gamma\leq \rho$ such that
$\max d^+(\cdots) \leq \rho$.
This will allow us to merge the universal block that should replace $\forall\mu$
directly into the already existing $\forall\rho\in\Q^{+*}$.

Now, for such a choice of $\gamma$,
and any $(w,b)\in \widetilde{W_\F^\epsilon}(\gamma)$,
there always exists some $\mu\in\widetilde{\M_\F^\B}(\epsilon)$
such that $d^+\left(\meas{(w,b)},\mu\right)\leq\rho$.
As $\pi_1^*\left(\mu\right)\in\M_\F^\B(\epsilon)$,
Equation~\ref{eqn:Stab2} applies to it,
and we can chose a corresponding pair
$\left(w_1,w_2\right)\in\left(\A^2\right)^{U_{\psi\left(\rho,\left|\A^2\right|,d\right)}}$.
Hence:
\[
\begin{array}{rcl}
d_{|\A|}^+\left(\meas{w_1},\meas{w}\right) &\leq&
d_{|\A|}^+\left(\meas{w_1},\pi_1^*\left(\mu\right)\right)
+d_{|\A|}^+\left(\pi_1^*\left(\mu\right),\meas{w}\right)\\
&\leq&d_{|\A|}^+\left(\meas{w_1},\pi_1^*\left(\mu\right)\right)+
d^+\left(\mu,\meas{(w,b)}\right)\\
&\leq &2\rho .
\end{array}
\]
Likewise, we have $\nu\in\M_\F$ such that $d_{|\A|}^+\left(\meas{w_2},\nu\right)\leq \rho$
in Equation~\ref{eqn:Stab2},
thus by the Covering Lemma~\ref{lem:CoveringMF} for $\M_\F$
we have a pattern $w_0\in W_\F(\rho)$ such $d_{|\A|}^+\left(\nu,\meas{w_0}\right)\leq \rho$,
hence $d_{|\A|}^+\left(\meas{w_2},\meas{w_0}\right)\leq 2\rho$.
The third inequality does not change, which concludes the implication
$(\ref{eqn:Stab2}\Rightarrow\ref{eqn:Stab3})$.

Now, suppose Equation~\ref{eqn:Stab3} is true and
let us prove $(\ref{eqn:Stab3}\Rightarrow\ref{eqn:Stab1})$.
Fix $\delta$, $\epsilon$ in the formula.
Consider any sequence $\rho_n\to 0$, and the corresponding $\gamma_n$ in the formula.

Let $\mu\in\M_\F^\B(\epsilon)$.
Using the Covering Lemma~\ref{lem:CoveringMFeps},
we know there exists a sequence 
$\left(w^n,b^n\right)\in\widetilde{\W_\F^\epsilon}\left(\gamma_n\right)$ such that
$d_{|\A|}^+\left(\mu,\meas{w^n}\right)\leq \gamma_n\leq \rho_n\to 0$.
At any rank, we may chose $w_0^n\in\W_\F\left(\rho_n\right)$
and $\left(w_1^n,w_2^n\right)\in\left(\A^2\right)^{
U_{\psi\left(\rho_n,\left|\A^2\right|,d\right)}}$ accordingly in Equation~\ref{eqn:Stab3}.
Up to extraction, $\meas{w_0^n}$ converges to $\nu\in\M_\F$
and $\meas{\left(w_1^n,w_2^n\right)}$ to $\lambda\in\M_{\A^2}$.

At the limit $\rho_n\to 0$,
the first inequality of Equation~\ref{eqn:Stab3} tells us that $\pi_1^*(\lambda)=\mu$.
Likewise, the second one tells us that $\pi_2^*(\lambda)=\nu\in\M_\F$,
thence $\lambda\in J\left(\mu,\M_\F\right)$.
The third inequality naturally becomes $\lambda(\Delta)\leq\delta$,
hence Equation~\ref{eqn:Stab1} holds true, the SFT is stable.
\end{proof}
\end{proposition}

\begin{theorem}[Upper Bound for Stability]
The problem $P_{stab}$ is in $\Pi_4$.

\begin{proof}
We just proved that $\F\in P_{halt}$ \emph{iff}
it satisfies Equation~\ref{eqn:Stab3}.
For the two blocks $d_{|\A|}^+\left(\meas{w_i},\meas{w}\right)\leq2\rho$
in Equation~\ref{eqn:Stab3},
we can replace $2\rho$ by $3\rho$ to have a strict inequality instead.
In particular, the proof of $(\ref{eqn:Stab3}\Rightarrow\ref{eqn:Stab1})$ applies
to this variant, so it is indeed an equivalent characterisation of stability.
The interest of this variant point of view is that,
as $d_{|\A|}^+\left(\meas{w_i},\meas{w}\right)$ is a computable real number,
$d_{|\A|}^+\left(\meas{w_i},\meas{w}\right)<3\rho$ becomes a \emph{semi-decidable} problem,
adding a countable existential quatifier that can be merged
into the $\exists\gamma$ block.

This variant formula starts with $\left[\forall\delta\in\Q^{+*},\exists\epsilon\in\Q^{+*},
\forall\rho\in\Q^{+*}, \exists\gamma\in\Q^{+*}\right]$,
\emph{i.e.}\ four alternating layers of countable quantifiers.
The following quantifiers are over finite computable sets,
and then the three inequalities are decidable.
Hence, this whole block can be decided in finite time.
\end{proof}
\end{theorem}

\newpage

\def\MR#1{\href{http://www.ams.org/mathscinet-getitem?mr=#1}{MR#1}}
\bibliographystyle{amsplain-nodash}
\bibliography{biblio.bib}

\end{document}